\documentclass[12pt,a4paper]{amsart}
\usepackage{amssymb,amscd}

\theoremstyle{plain}
\newtheorem{theorem}{Theorem}[subsection]
\newtheorem{lemma}[theorem]{Lemma}
\newtheorem{proposition}[theorem]{Proposition}
\newtheorem{corollary}[theorem]{Corollary}

\theoremstyle{definition}

\newtheorem{remark}[theorem]{Remark}
\newtheorem{remarks}[theorem]{Remarks}
\newtheorem{example}[theorem]{Example}

\numberwithin{equation}{subsection}

\newcommand\bA{{\mathbb A}}

\newcommand\bF{{\mathbb F}}
\newcommand\bG{{\mathbb G}}
\newcommand\bI{{\mathbb I}}
\newcommand\bP{{\mathbb P}}

\newcommand\cA{{\mathcal A}}

\newcommand\cE{{\mathcal E}}
\newcommand\cF{{\mathcal F}}

\newcommand\cL{{\mathcal L}}
\newcommand\cN{{\mathcal N}}
\newcommand\cO{{\mathcal O}}
\newcommand\cQ{{\mathcal Q}}
\newcommand\cS{{\mathcal S}}
\newcommand\cT{{\mathcal T}}

\newcommand\fc{{\mathfrak c}}

\newcommand\fg{{\mathfrak g}}
\newcommand\fh{{\mathfrak h}}

\newcommand\fm{{\mathfrak m}}
\newcommand\fn{{\mathfrak n}}

\newcommand\aff{\operatorname{aff}}

\newcommand\diag{\operatorname{diag}}

\newcommand\op{\operatorname{op}}

\newcommand\red{\operatorname{red}}

\newcommand\Aut{\operatorname{Aut}}

\newcommand\GL{\operatorname{GL}}
\newcommand\Gr{\operatorname{Gr}}

\newcommand\PGL{\operatorname{PGL}}

\title{Log homogeneous varieties}

\author{Michel~Brion}
\address{Universit\'e de Grenoble I\\
D\'epartement de Math\'ematiques\\
Institut Fourier, UMR 5582 du CNRS\\
38402 Saint-Martin d'H\`eres Cedex, France}
\email{Michel.Brion@ujf-grenoble.fr}

\begin{document}
 
\begin{abstract}
Given a complete nonsingular algebraic variety $X$ and a divisor $D$
with normal crossings, we say that $X$ is log homogeneous with
boundary $D$ if the logarithmic tangent bundle 
$T_X( - \log D)$ is generated by its global sections. Then the
Albanese morphism $\alpha$ turns out to be a fibration with fibers being 
spherical (in particular, rational) varieties. It follows that all
irreducible components of $D$ are nonsingular, and any partial
intersection of them is irreducible. Also, the image of $X$ under the
morphism $\sigma$ associated with $- K_X - D$ is a spherical
variety, and the irreducible components of all fibers of $\sigma$ are
equivariant compactifications of semiabelian varieties. Generalizing
the Borel--Remmert structure theorem for homogeneous varieties, we
show that the product morphism $\alpha \times \sigma$ is surjective,
and the irreducible components of its fibers are toric varieties. We
reduce the classification of log homogeneous varieties to a problem
concerning automorphism groups of spherical varieties, that we solve
under an additional assumption.
\end{abstract}

\maketitle

\tableofcontents

\setcounter{section}{-1}

\section{Introduction}
\label{sec:introduction}

Consider a complete nonsingular algebraic variety $X$ over the field
of complex numbers. If $X$ is homogeneous, then its structure is well
understood by a classical result of Borel and Remmert: $X$ is the
product of an abelian variety (the Albanese variety) and a flag
variety (the basis of the ``Tits fibration''). 

\medskip

More generally, one would like to classify the almost homogeneous
varieties $X$, i.e., those where a connected algebraic group $G$ acts
with an open orbit $X_0$. However, this central problem of equivariant
geometry seems to be too general, and complete results are only known
under assumptions of smallness (in various senses) for the boundary 
$X \setminus X_0$. We refer to \cite{HO84} for an exposition of such
classification results, in the setting of holomorphic actions on
complex analytic spaces.

\medskip

In another direction, there is a well-developed structure theory for
certain classes of almost homogeneous varieties, where the acting
group and the open orbit are prescribed: torus embeddings, spherical
embeddings and, more generally, equivariant embeddings of a
homogeneous space $X_0$ under a reductive group $G$ (see the
recent survey \cite{Ti06}). A subclass of special interest is that of
$G$-regular varieties in the sense of \cite{BDP90}; in loose words, 
their orbit structure is that of a nonsingular toric variety. Another
important subclass consists of semiabelic varieties; in loose words
again, they are toric bundles over abelian varieties. (Semiabelic
varieties are introduced in \cite{Al02}, where their role in
degenerations of abelian varieties is investigated.)

\medskip

These restrictive assumptions on $G$ and $X_0$ are convenient in the
setting of algebraic transformation groups, but somewhat unnatural 
from a geometric viewpoint. This motivates the search for a class of
almost homogeneous varieties containing all regular varieties
under a reductive group and all semiabelic varieties, having geometric
significance and an accessible structure. Also, it is natural to
impose that the boundary be a divisor with normal crossings.

\medskip

In the present paper, we introduce the class of log homogeneous
varieties and show that it fulfills the above requirements. 
Given a normal crossings divisor $D$ in $X$, we say that $X$ 
is log homogeneous with boundary $D$ if the associated logarithmic
tangent bundle $T_X (- \log D)$ is generated by its global
sections. It follows readily that $X$ is almost homogeneous under the
connected automorphism group $G := \Aut^0(X,D)$, with boundary being
$D$. More generally, the $G$-orbits in $X$ are exactly the strata
defined by $D$ (Prop.~\ref{prop:crit}); in particular, their
number is finite. 

\medskip

All $G$-regular varieties are log homogeneous, as shown in \cite{BB96}. 
Semiabelic varieties satisfy a stronger property, namely, the bundle
$T_X(- \log D)$ is trivial; we then say that $X$ is log
parallelizable with boundary $D$. Conversely, log parallelizable
varieties are semiabelic by a theorem of Winkelmann \cite{Wi04}
which has been the starting point for our investigations. The class of
log homogeneous varieties turns out to be stable under several natural
operations: induction (Prop.~\ref{prop:Gf}), equivariant blow-ups
(Prop.~\ref{prop:mor}), \'etale covers (Prop.~\ref{prop:ec})
and taking invariant subvarieties or irreducible components of fibers
of morphisms (Cor.~\ref{cor:lhr}). 

\medskip

Any log homogeneous variety $X$ comes with two natural morphisms
which turn out to play opposite roles. The first one is the Albanese
map $\alpha : X \to \cA(X)$; in our algebraic setting, this is the 
universal map to an abelian variety. We show that $\alpha$ is a  
fibration with fibers being log homogeneous under the
maximal connected affine subgroup $G_{\aff}$ of $G$, and spherical
under any Levi subgroup of $G_{\aff}$; this yields our main structure
theorem (Thm.~\ref{thm:str}). 

\medskip

The second morphism is the Tits map $\tau$ to a Grassmannian, defined
by the global sections of the logarithmic tangent bundle. 
Let $\sigma : X \to \cS(X)$ be the Stein factorization of $\tau$, so
that $\sigma$ is the morphism associated with the globally generated
divisor $- K_X - D$ (the determinant of $T_X( - \log D)$). We show
that $\cS(X)$ is also a spherical variety under any Levi subgroup of
$G_{\aff}$, and the irreducible components of fibers of $\sigma$ are
semiabelic varieties (Prop.~\ref{prop:sf}). It follows that the
product morphism $\alpha \times \sigma$ is surjective with fibers
being finite unions of toric varieties (Thm.~\ref{thm:atm}). This
generalizes the Borel--Remmert theorem; note, however, that $\cS(X)$
may be singular, see Remark \ref{rem:fin}(ii).

\medskip

Actually, these results hold in a slightly more general setting,
namely, for a faithful action of a connected algebraic group $G$ on
$X$ preserving a normal crossings divisor $D$ and such that the
associated global vector fields generate $T_X (- \log D)$. Then $G$
has the same orbits as the full group $\Aut^0(X,D)$, but may well be
strictly contained in that group. For example, let $X$ be the
projective $n$-space and $D$ the union of $m$ hyperplanes in general
position. If $m \le n$, then $X$ is log homogeneous with boundary $D$
under the action of $\Aut^0(X,D)$. If, in addition, $n - m$ is odd and
at least $3$, then $X$ is still log homogeneous with boundary $D$ for
a smaller group $G$, namely, the largest subgroup of $\Aut^0(X,D)$ that
acts on the intersection of the $m$ hyperplanes as the projective
symplectic group.

\medskip

This article is organized as follows. Section 1 contains a number of
useful preliminary results on algebraic groups and their homogeneous
spaces, that we did not see explicitly stated in the literature. We
discuss algebraic analogues of the Albanese and Tits fibrations, which
are classical tools in the theory of complex homogeneous manifolds. In
our algebraic setting, the Borel--Remmert theorem admits a very short
proof (Theorem \ref{thm:br}).

\medskip

In Section 2, we consider nonsingular almost homogeneous
varieties, possibly not complete. We obtain a criterion for
log homogeneity resp. parallelizability, as well as stability
properties and simple local models. We also show that the Albanese map
is a fibration; this yields an important reduction to the case where
the acting group is affine. Finally, we obtain an algebraic version
of Winkelmann's theorem, again with a very short proof (Theorem
\ref{thm:lpt}).

\medskip

The final Section 3 is devoted to complete log homogeneous
varieties. After some preliminary lemmas about their relations to
spherical varieties, we describe their structure and analyze the Tits
map and its Stein factorization $\sigma$. We also show that the group
of equivariant automorphisms is an extension of the Albanese variety
by a diagonalizable group, so that its connected component is 
a semiabelian variety, acting on the general fibers of $\sigma$ with
an open orbit (Props.~\ref{prop:ea} and \ref{prop:cc}). In fact, the
relative logarithmic tangent bundle of $\sigma$ is trivial, with fiber
being the Lie algebra of the equivariant automorphism group
(Prop.~\ref{prop:sf}). Finally, we introduce the subclass of
strongly log homogeneous varieties, closely related to regular
varieties under a reductive group. We obtain a simple characterization
of this subclass (Prop.~\ref{prop:slh}) and we apply it to the
classification of complete log homogeneous surfaces
(Prop.~\ref{prop:cs}).

\medskip

Our arguments are purely algebraic; they involve results of the theory
of algebraic transformation groups, and some basic notions borrowed
from logarithmic birational geometry (for which we refer to
\cite{Ma02}) and holomorphic transformation groups (see
\cite{Ak95,HO84}). Our results represent only the first steps in the
understanding of log homogeneous varieties. They raise many open
problems, among which the most important ones seem to be:

\smallskip

\noindent
A) Characterize those homogeneous spaces that admit a log homogeneous
completion. 

\smallskip

\noindent
B) Classify the triples $(L,Y,N)$, where $L$ is a connected reductive
group, $Y$ is a complete nonsingular spherical $L$-variety, and $N$ is
a unipotent group of automorphisms of $Y$, normalized by $L$ and such
that $Y$ is log homogeneous under the semidirect product $LN$.

\medskip

By Theorem \ref{thm:str}, Problem A may be reduced to homogeneous
spaces $X_0$ under an affine group $G$, and then a necessary condition
is that $X_0$ be spherical under a Levi subgroup of $G$. We do not
know if this condition is not sufficient.

\medskip

Problem B is equivalent to the classification of complete log
homogeneous varieties, by Theorem \ref{thm:str} again. It fits into
the problem of describing automorphism groups of complete nonsingular
spherical varieties, which is very much open in general (see
\cite[Sec.~3.4]{Od88} for the toric case, and \cite{Br06} for the
regular case).

\medskip

Another natural question concerns the case where the field of complex
numbers is replaced with an algebraically closed field of positive
characteristics. Here a number of results carry over with only small
changes, but the Borel--Remmert theorem requires strong additional
assumptions of separability. Also, some important ingredients of
Section 3 (the existence of a Levi decomposition, the local structure
of spherical varieties, Knop's vanishing theorem for the logarithmic
tangent sheaf) are not available in this setting.

\bigskip

\noindent
{\bf Notation and conventions.} 
Throughout this article, we consider algebraic varieties and algebraic
groups over an algebraically closed field $k$ of characteristic
zero. By a variety, we mean an integral separated scheme of finite
type over $k$, and by an algebraic group, a group scheme of finite
type over $k$. Morphisms are understood to be morphisms of varieties
over $k$; points are understood to be $k$-rational. As a general
reference for algebraic geometry, we use the book \cite{Ha77}, and
\cite{DG70} for algebraic groups.

Any algebraic group $G$ is nonsingular, since $k$ has char.~$0$.
Thus, the neutral component $G^0$, that is, the connected component of
$G$ containing the identity element, is a nonsingular variety. Also,
recall that being affine or linear are equivalent properties of $G$
(or $G^0$). But we shall consider algebraic groups that are not
necessarily affine, e.g., abelian varieties for which we refer to the
expository article \cite{Mi86}.

\section{Homogeneous varieties}
\label{sec:hs}

\subsection{Algebraic groups}
\label{subsec:ag}

Let $G$ be a connected algebraic group. By a result of Chevalley (see
\cite{Co02} for a modern proof), there exists a unique closed
connected normal affine subgroup $G_{\aff} \subseteq G$ such that the
quotient group $G/G_{\aff}$ is an abelian variety. We denote by
$$
\alpha: G \to G/G_{\aff} =: \cA(G)
$$
the quotient morphism.

Recall also that any morphism from a connected affine algebraic group
to an abelian variety is constant, see \cite[Cor.~3.9]{Mi86}. It follows
that $G_{\aff}$ contains every closed connected affine subgroup of
$G$. Moreover, every morphism (of varieties) $f : G \to A$, where $A$
is an abelian variety, factors uniquely as $\varphi \circ \alpha$,
where $\varphi : \cA(G) \to A$ is a morphism. By \cite[Cor.~2.2]{Mi86},
$\varphi$ is the composition of a translation and a group
homomorphism. So $\alpha$ is the {\it Albanese morphism} of the
variety $G$, as defined in \cite{Se58}. In particular, $\cA(G)$
depends only on $G$ regarded as a variety.

We record the following easy result, where we denote the center
of $G$ by $C(G)$, or simply by $C$ if this yields no confusion. 
Also, recall that an {\it isogeny} is a surjective homomorphism of
algebraic groups, with finite kernel.

\begin{lemma}\label{lem:G} 
The algebraic groups $G/C$ and $G/C^0$ are both affine, and  
\begin{equation}\label{eqn:gz}
G = G_{\aff} C^0.
\end{equation}
Moreover, $C^0_{\aff}$ is the neutral component of both groups
$G_{\aff} \cap C$ and $G_{\aff} \cap C^0$, and the natural map 
$\cA(C^0) \to \cA(G)$ is an isogeny.
\end{lemma}

\begin{proof}
Via the adjoint representation of $G$ in its Lie algebra $\fg$, the
quotient group $G/C$ is isomorphic to a closed subgroup of
$\GL(\fg)$. In particular, $G/C$ is affine. 

The natural map $G/C^0 \to G/C$ is the quotient by the finite group
$C/C^0$, and hence is an affine morphism. Thus, $G/C^0$ is affine as
well.

Note that $G_{\aff} C^0$ is a closed normal subgroup of $G$, and the
corresponding quotient 
$$
G/G_{\aff} C^0 \cong (G/C^0)/(G_{\aff} C^0/C^0)
$$
is affine. On the other hand, $G/G_{\aff} C^0$ is a quotient of the 
abelian variety $\cA(G)$, and hence is complete. This yields
(\ref{eqn:gz}). 

The remaining assertions follow readily from the definitions.
\end{proof}

\subsection{The Albanese fibration}
\label{subsec:af1}

Let $X$ be a $G$-variety and $X_0 \subseteq X$ a $G$-orbit. 
Choosing a base point $x \in X_0$, we identify the homogeneous variety
$X_0 = G \cdot x$ with the homogeneous space $G/H$, where $H$ denotes
the isotropy group $G_x$.

Since $G/G_{\aff}$ is an abelian variety, the product 
$$
I(X_0) := G_{\aff} H
$$ 
is a closed normal subgroup of $G$, independent of the choice of the
base point $x$, and the quotient $G/I(X_0)$ is an abelian variety as
well. Moreover, the natural map 
$$
\alpha : X_0 = G/H \to G/G_{\aff} H = X_0/G_{\aff} =: \cA(X_0)
$$
is the Albanese morphism of $X_0$. Note that $\alpha$ is a
$G$-equivariant fibration with fiber
$$
G_{\aff} H/H \cong G_{\aff}/(G_{\aff}\cap H),
$$
and this $G_{\aff}$-homogeneous variety is also independent of the
choice of $x$.

Further properties of $H$ and $I(X_0)$ are gathered in the following:

\begin{lemma}\label{lem:H}
Let $X$ be a variety on which $G$ acts faithfully, and 
$X_0 = G \cdot x \subseteq X$ a $G$-orbit. Then the isotropy group 
$H = G_x$ is affine; equivalently,  
$H^0 \subseteq G_{\aff}$. Moreover, there exists a finite subgroup 
$F \subseteq C^0$ such that 
\begin{equation}\label{eqn:if}
I(X_0) = G_{\aff} F.
\end{equation} 
\end{lemma}

\begin{proof}
Note that $H$ acts on the local ring $\cO_x$ and on its
quotients $\cO_x/\fm_x^n$ by powers of the maximal ideal; these
quotients are finite-dimensional $k$-vector spaces. Let $H_n$ be the
kernel of the resulting homomorphism $H \to \GL(\cO_x/\fm_x^n)$. Then
the $H_n$ form a decreasing sequence of closed subgroups of $H$, and
their intersection acts trivially on $\cO_x$. Since the $G$-action is 
faithful, this intersection is trivial, and hence so is $H_n$ for 
$n \gg 0$. Thus, the algebraic group $H$ is affine.

As a consequence, $I(X_0)^0 = G_{\aff}$, and $I(X_0)$ is affine as
well. Moreover,
$$
I(X_0) = G_{\aff} (C^0 \cap I(X_0))
$$
by (\ref{eqn:gz}). Thus, $C^0 \cap I(X_0)$ is a commutative
affine algebraic group: it is the direct product of a connected
algebraic group (contained in $G_{\aff}$) with a finite group.
\end{proof}

\subsection{The Tits fibration}
\label{subsec:tf}

Consider again a $G$-variety $X$. The tangent sheaf 
$$
\cT_X = Der_k(\cO_X)
$$ 
is a $G$-linearized sheaf on $X$, equipped with a $G$-equivariant map
\begin{equation}\label{eqn:OP}
\op_X: \fg \to \Gamma(X,\cT_X).
\end{equation}
This yields a morphism of $G$-linearized sheaves
\begin{equation}\label{eqn:op}
op_X: \cO_X \otimes \fg \to \cT_X.
\end{equation}
Clearly, $op_X$ is surjective if and only if $X$ is
$G$-homogeneous. Under this assumption, the kernel of $op_X$ is a
$G$-linearized locally free sheaf, and its fiber at any point $x$ is
the isotropy Lie subalgebra $\fg_x \subseteq \fg$. Likewise, the fiber
at $x$ of the locally free $G$-linearized sheaf $\cT_X$ is the
quotient $\fg/\fg_x$. In particular, $op_X$ is an isomorphism if and
only if $X$ is the quotient of $G$ by a finite subgroup; we then say
that $X$ is $G$-{\it  parallelizable}.

If $X$ is homogeneous under a faithful action of $G$, then $\fg_x$ is
contained in the Lie algebra $\fg_{\aff}$ of $G_{\aff}$ by Lemma
\ref{lem:H}. So we obtain a morphism 
$$
\tau: X \to \cL, \quad x \mapsto \fg_x \; ,
$$
where $\cL$ denotes the scheme of Lie subalgebras of $\fg_{\aff}$
(alternatively, we may consider the Grassmann variety of subspaces of
$\fg_{\aff}$). Clearly, $\tau$ is $G$-equivariant, where $G$ acts on
$\cL$ via its adjoint action on $\fg_{\aff}$; in particular, $\tau$
is invariant under $C$. Thus, the image of $\tau$ is a unique
$G_{\aff}$-orbit that we denote by $\cL(X)$, and $\tau : X \to \cL(X)$
is a $G$-equivariant fibration, called the {\it Tits fibration}.

Choosing a base point $x \in X$ with isotropy group $H$, we obtain
with obvious notation
$$
\cL(X) = G \cdot \fh \cong G/N_G(\fh) = G/N_G(H^0) \cong
G_{\aff}/N_{G_{\aff}}(H^0) \cong G_{\aff}/N_{G_{\aff}}(\fh),
$$
where the second isomorphism follows from Lemma \ref{lem:G}. Thus,
the fiber of $\tau$ at $x$ is the homogeneous space $N_G(\fh)/H$, 
quotient of the (not necessarily connected) algebraic group
$N_G(\fh)/H^0$ by the finite subgroup $H/H^0$. In particular, all
connected components of fibers are parallelizable varieties. 

\begin{remarks}\label{rem:tf}
(i) The Tits fibration may be defined, more generally, for a
homogeneous $G$-variety $X$ on which the action is not necessarily
faithful: just replace $\fg_x$ with $\fg_{\aff,x}$ (the isotropy Lie
algebra of $x$ in $\fg_{\aff}$). But several statements take a much
simpler form in the setting of faithful actions.

\smallskip

\noindent
(ii) Unlike the Albanese fibration, the Tits fibration may depend on
the $G$-action on $X$. Consider, for example, a connected affine group
$G$ regarded as a homogeneous $G$-variety under its action via left
multiplication; then the Tits fibration is trivial. But this fibration
is nontrivial if $G$ is regarded as a homogeneous 
$G \times G$-variety under its action via left and right multiplication,
and $G$ is noncommutative.

\smallskip

\noindent
(iii) The Tits fibration yields an exact sequence of $G$-linearized
locally free sheaves on $G/H$
$$
0 \to \cT_{\tau} \to \cT_X \to \tau^*\cT_{\cL(X)} \to 0,
$$
where $\cT_{\tau}$ denotes the relative tangent sheaf. By taking fibers
at the base point, this corresponds to an exact sequence of
$H$-modules
$$
0 \to \fn_{\fg}(\fh)/\fh \to \fg/\fh \to \fg/\fn_{\fg}(\fh) \to 0.
$$
Note that $H$ acts on $\fn_{\fg}(\fh)/\fh)$ via its finite quotient
$H/H^0$. Thus, the pull-back of $T_{\tau}$ under the finite cover
$G/H^0 \to G/H$ is trivial.

\smallskip

\noindent
(iv) Likewise, the Albanese fibration yields an exact sequence
of $G$-linearized locally free sheaves
$$
0 \to \cT_{\alpha} \to \cT_X \to \alpha^*\cT_{\cA(X)} \to 0
$$
corresponding to the exact sequence of $H$-modules
$$
0 \to \fg_{\aff}/\fh \to \fg/\fh \to \fg/\fg_{\aff} \to 0.
$$
Note that the $H$-module $\fg/\fg_{\aff}$ is trivial; this also
follows from the triviality of the tangent sheaf of the abelian
variety $\cA(X)$. 
\end{remarks}

\subsection{The product fibration}
\label{subsec:pf}

We now consider the product map 
$$
\alpha \times \tau : X \to \cA(X) \times \cL(X),
$$
where $X \cong G/H$.

\begin{lemma}\label{lem:at}
With the above notation, $\alpha \times \tau$
is surjective; in particular, it is a $G$-equivariant fibration. Its
fiber at the base point of $G/H$ is the homogeneous space
$N_{G_{\aff}}(\fh)/(G_{\aff}\cap H)$, quotient of the affine
algebraic group $N_{G_{\aff}}(\fh)/H^0$ by the finite subgroup
$(G_{\aff}\cap H)/H^0$.

In particular, the morphism $\alpha \times \tau$ is affine, and the 
connected components of its fibers are parallelizable varieties.
\end{lemma}

\begin{proof}
Any fiber of $\alpha$ is a unique $G_{\aff}$-orbit, and hence is
mapped onto $\cL(X)$ by $\tau$. Thus, $\alpha \times \tau$ is
surjective. In particular, $\cA(X) \times \cL(X)$ is a unique
$G$-orbit. Since $\cA(X) = G/G_{\aff}H$ and $\cL(X) = G/N_G(\fh)$,
it follows that 
$$
\cA(X) \times \cL(X) = G/(G_{\aff}H \cap N_G(\fh))
= G/N_{G_{\aff}}(\fh) H.
$$
Thus, the fiber of $\alpha \times \tau$ equals
$N_{G_{\aff}}(\fh) H/H \cong N_{G_{\aff}}(\fh)/(G_{\aff}\cap H)$.
\end{proof}

\begin{remark}
As in Remark \ref{rem:tf}(iii), the relative tangent sheaf of
$\alpha \times \tau$ corresponds to the $H/H^0$-module 
$\fn_{\fg_{\aff}}(\fh)/\fh$. Thus, the pull-back of this sheaf under
the finite cover $G/H^0 \to G/H$ is trivial.
\end{remark}

The above lemma yields an algebraic version of the Borel--Remmert
structure theorem for compact homogeneous K\"ahler manifolds, see
e.g.~\cite[\S 3.9]{Ak95}:   

\begin{theorem}\label{thm:br}
For a complete homogeneous $G$-variety $X$, the map 
$\alpha \times \tau$ is an isomorphism.
\end{theorem}

\begin{proof}
We may assume that $G$ acts faithfully on $X = G/H$. Since 
the fiber of $\alpha$ is $G_{\aff}/(G_{\aff}\cap H)$, then 
$G_{\aff}\cap H$ is a parabolic subgroup of $G_{\aff}$. In particular,
$G_{\aff}\cap H$ is connected and equals its normalizer in
$G_{\aff}$. It follows that 
$$
H^0 = G_{\aff} \cap H = N_{G_{\aff}}(H^0) = N_{G_{\aff}}(\fh)
$$
so that the morphism $\alpha \times \tau$ is bijective by Lemma
\ref{lem:at}. Since $\cA(X) \times \cL(X)$ is nonsingular, 
$\alpha \times \tau$ is an isomorphism.
\end{proof}

\section{Log homogeneous varieties}
\label{sec:pairs}

\subsection{Definitions and basic properties}
\label{subsec:dbp}
 
Consider a pair $(X,D)$, where $X$ is a nonsingular variety and
$D \subset X$ is a {\it divisor with normal crossings}. That is, $D$
is an effective divisor and its local equation at an arbitrary point
$x \in X$ decomposes in the completed local ring $\widehat{\cO}_x$
into a product $t_1 \cdots t_r$, where $t_1,\ldots,t_r$ form part of a
regular system of parameters $(t_1,\ldots,t_n)$ of $\widehat{\cO}_x$. 
Let  
$$
\cT_X (- \log D) \subseteq \cT_X = Der_k(\cO_X)
$$ 
be the subsheaf consisting of those derivations that
preserve the ideal sheaf $\cO_X(-D)$. One easily checks that the 
{\it logarithmic tangent sheaf} $\cT_X(-\log D)$ is a locally free
sheaf of Lie subalgebras of $\cT_X$, having the same restriction to $X
\setminus D$, and hence the same rank $n = \dim(X)$.
The dual of $\cT_X(-\log D)$ is the sheaf $\Omega^1_X( \log D)$ of
logarithmic differential forms, that is, of differential $1$-forms on
$X \setminus D$ having at most simple poles along $D$. The top
exterior power $\wedge^n \cT_X(-\log D)$ is the invertible sheaf
$\cO_X(-K_X - D)$, where $K_X$ denotes the canonical divisor.

If $D$ is defined at $x$ by the equation $t_1\cdots t_r=0$ as above,
then a local basis of $\cT_X(- \log D)$ (after localization and
completion at $x$) consists of
$$
t_1 \partial_1, \ldots, t_r \partial_r, \partial_{r+1}, \ldots, 
\partial_n,
$$ 
where $(\partial_1, \ldots, \partial_n)$ is the local basis of $\cT_X$
dual to the local basis $(dt_1, \ldots, dt_n)$ of $\Omega^1_X$.

Next, assume that a connected algebraic group $G$ acts on $X$ and
preserves $D$; we say that $(X,D)$ is a $G$-{\it pair}. Then the 
sheaf $\cT_X( -\log D)$ is $G$-linearized and the map 
$\op_X$  of (\ref{eqn:OP}) factors through a map  
\begin{equation}\label{eqn:OPD}
\op_{X,D} : \fg \to \Gamma(X,\cT_X(- \log D)).
\end{equation}
This defines a morphism of $G$-linearized sheaves
\begin{equation}\label{eqn:opd}
op_{X,D}: \cO_X \otimes \fg \to \cT_X(-\log D).
\end{equation}
We say that the pair $(X,D)$ is 
{\it homogeneous ({\rm resp.}~parallelizable) under} $G$, if
$op_{X,D}$ is surjective (resp.~an isomorphism).

More generally, given a locally closed $G$-stable subvariety 
$S \subseteq X$, we say that $(X,D)$ is {\it homogeneous 
({\rm resp.}~parallelizable) under $G$ along $S$}, 
if $op_{X,D}$ is surjective (resp.~an isomorphism) at all points of $S$.

\begin{remarks}\label{rem:lh}
(i) For any homogeneous $G$-pair $(X,D)$, the complement 
$X_0 := X \setminus D$ is a unique $G$-orbit, since $op_{X_0}$ is
surjective. In particular, the $G$-variety $X$ is almost
homogeneous, and $D$ is uniquely determined by the pair $(X,G)$. We
then say that $X$ is a 
{\it log homogeneous $G$-variety with boundary $D$}. Log
parallelizable $G$-varieties are defined similarly.

\smallskip

\noindent
(ii) If $X$ is complete, then the connected automorphism group
$\Aut^0(X)$ is algebraic with Lie algebra being $\Gamma(X,\cT_X)$, see
e.g.~\cite{MO67}. It follows that $\Aut^0(X,D)$ is algebraic as well,
with Lie algebra being $\Gamma(X,\cT_X(- \log D))$. Thus, $(X,D)$ is 
homogeneous (resp.~parallelizable) under some group $G$ if and
only if the sheaf $\cT_X(- \log D)$ is generated by its global
sections (resp.~is trivial). We then say that the pair $(X,D)$ is
homogeneous (resp.~parallelizable) without specifying the group $G$.

\smallskip

\noindent
(iii) The log homogeneous $G$-pairs $(X,D)$, where $X$ is a curve and
$G$ acts faithfully, are easily classified. One obtains the following
list of triples $(X,D,G)$:

\smallskip

\noindent
(parallelizable) $(G,\emptyset,G)$, where $G$ is a connected algebraic
group of dimension $1$, i.e., the additive group $\bG_a$, the
multiplicative group $\bG_m$, or an elliptic curve.

\smallskip

\noindent
(homogeneous) $(\bP^1,\emptyset,\PGL_2)$ and $(\bA^1,\emptyset,B)$,
where $B$ denotes the automorphism group of the affine line, i.e., the
Borel subgroup of $\PGL_2$ that fixes $\infty \in \bP^1$.

\smallskip

\noindent
(log parallelizable) $(\bA^1,\{0\},\bG_m)$ and $(\bP^1,\{0,\infty\},\bG_m)$.

\smallskip

\noindent
(log homogeneous) $(\bP^1,\{\infty\},B)$.

\smallskip

The complete log homogeneous surfaces will be listed in  
Proposition \ref{prop:cs}.
\end{remarks}

We shall obtain a criterion for an arbitrary $G$-pair $(X,D)$ to be
homogeneous or parallelizable; to state it, we first recall the
definition of a stratification of $X$ associated with $D$. Let 
$X^0 :=X$, $X^1 := D$ and define inductively $X^r$ to be the singular 
locus of $X^{r-1}$ for $r \ge 2$. Then the strata are the connected
components of the locally closed subvarieties $X^r \setminus X^{r+1}$,
where $r = 0,1,\ldots$  Every stratum is nonsingular and $G$-stable;
the open stratum is $X \setminus D$.

Next, we analyze the normal space to a stratum $S$ at a point $x$. 
Let again $t_1 \cdots t_r=0$ be a local equation of $D$ at $x$. 
Then the ideal sheaf of $S$ in $X$ is generated (after
localization and completion at $x$) by $t_1, \ldots, t_r$; in other
words, $S$ is locally the complete intersection of all branches of $D$. 
As a consequence, the normal space $N_{S/X,x}$ admits a canonical
splitting into a direct sum of lines $L_1, \ldots, L_r$; the
corresponding hyperplanes $\bigoplus_{j\ne i} L_j$ are the normal
spaces to these branches. The isotropy group $G_x$ preserves $D$ and
permutes its branches at $x$, and hence the lines $L_1, \ldots, L_r$;
the connected component $G_x^0$ stabilizes each line. Thus, the
representation of $G_x$ in $N_{S/X,x}$ yields a homomorphism
\begin{equation}\label{eqn:rho}
\rho_x: G_x^0 \to \GL(L_1) \times \cdots \times \GL(L_r) 
\cong \bG_m^r \; ,
\end{equation}
where $\bG_m^r$ denotes the product of $r$ copies of the
multiplicative group.

We may now state our criterion:

\begin{proposition}\label{prop:crit}
The following conditions are equivalent for a $G$-pair $(X,D)$ and a
stratum $S$:

\smallskip

\noindent
{\rm (i)} $(X,D)$ is homogeneous (resp.~parallelizable) under
$G$ along $S$.

\smallskip

\noindent
{\rm (ii)} $S$ is a unique $G$-orbit and for any $x \in S$, the
homomorphism $\rho_x$ is surjective (resp.~an isogeny).

\smallskip

If one of these conditions holds, then the sequence
\begin{equation}\label{eqn:ex}
\CD
0 @>>> \fg_{(x)} @>>> \fg @>{\op_{x,X,D}}>> T_x X(- \log D) @>>> 0
\endCD
\end{equation}
is exact for any $x \in S$, where $\fg_{(x)}$ denotes the kernel of
the representation of the isotropy Lie algebra $\fg_x$ in the normal
space $N_{S/X,x}$ and $T_x X(- \log D)$ denotes the fiber of
$\cT_X(-\log D)$ at $x$.
\end{proposition}

\begin{proof}
Clearly, $\cT_X(-\log D)$ preserves the ideal sheaf of $S$ in $X$.
This yields a morphism 
$\cT_X(-\log D)\vert_S \to \cT_S$ 
and, in turn, a linear map 
$$
p : T_x X(-\log D) \to T_x S
$$ 
between fibers. Given a regular system of parameters $(t_1,\ldots,t_n)$
as above, the map $p$ is just the projection of the space
$k t_1 \partial_1 \oplus \cdots \oplus k t_r \partial _r
\oplus k \partial_{r+1} \oplus \cdots \oplus k \partial_n$
onto its subspace 
$k \partial_{r+1} \oplus \cdots \oplus k \partial_n$. 
Thus, $p$ fits into an exact sequence
$$
0 \to \bigoplus_{i=1}^r k t_i \partial_i \to T_x X(-\log D) \to T_x S \to 0.
$$
Moreover, composing the map $\op_{x,X,D}: \fg \to T_x X(-\log D)$ with
$p$ yields the map $\op_{x,S} : \fg \to T_x S$ which factors through
an injective map
$$
i:\fg/\fg_x \to T_x S.
$$ 
So we obtain a commutative diagram of exact sequences
$$
\CD
0 @>>> \fg_x @>>> \fg           @>>> \fg/\fg_x @>>> 0 \\
  & &            @VVV               @VVV           @V{i}VV   \\
0 @>>> \bigoplus_{i=1}^r k t_i \partial_i   @>>> T_x X(-\log D) 
@>{p}>> T_x S    @>>> 0, \\
\endCD 
$$
where the left vertical map may be identified with the differential
$d\rho_x : \fg_x \to k^r$. As a consequence, the surjectivity of the
middle vertical map is equivalent to the surjectivity of both maps
$d\rho_x$ and $\op_{x,S}$. If the latter condition holds, then we
obtain an exact sequence
$$
0 \to \ker(d\rho_x) \to \fg \to T_x X(-\log D) \to 0.
$$
This implies all assertions.
\end{proof}

A direct consequence of this criterion is the following:

\begin{corollary}\label{cor:fo}
Let $(X,D)$ be a homogeneous pair under $G$ acting faithfully. Then:

\smallskip

\noindent
{\rm (i)} The $G$-orbits in $X$ are exactly the strata; in particular,
their number is finite.

\smallskip

\noindent
{\rm (ii)} The Tits fibration of the open orbit $X \setminus D$
extends to a $G$-equivariant morphism
$$
X \to \cL, \quad x \mapsto \fg_{(x)}
$$
which is invariant under $C$.
\end{corollary}

Next, we compare our notion of log homogeneity to the (earlier) notions
of pseudofreeness and regularity. Recall from \cite[Sec.~2]{Kn94b}
that a nonsingular $G$-variety $X$ is {\it pseudofree} if the image of 
$op_X$ (i.e., the coherent subsheaf of $\cT_X$ generated by the image
of $\fg$) is locally free. Clearly, any log homogeneous variety is
pseudofree, but the converse does not hold in general. For example,
the projective line where the additive group acts by affine
translations is pseudofree but not log homogeneous.

Also, recall from \cite{BDP90} that our nonsingular $G$-variety $X$ is 
$G$-{\it regular} if it satisfies the following three conditions:

\smallskip

\noindent
(i) $X$ contains an open $G$-orbit whose complement is a union of
irreducible divisors with normal crossings (the {\it boundary divisors}).

\smallskip

\noindent
(ii) Any $G$-orbit closure in $X$ is the transversal intersection of
those boundary divisors that contain it.

\smallskip

\noindent
(iii) For any $x \in X$, the normal space $N_{G \cdot x/X,x}$ contains
an open orbit of $G_x$.

\begin{corollary}\label{cor:rlh}
Every regular $G$-variety is log homogeneous.
\end{corollary}

\begin{proof}
By (i), the complement of the open orbit is a divisor with normal
crossings. Moreover, by (ii), the strata are exactly the $G$-orbits,
and each normal space $N_{G \cdot x/X,x} \cong \bA^r$ is a direct sum
of lines $L_1,\ldots,L_r$, where each $L_i$ is preserved under $G_x$
(as $L_i$ is the intersection of normal spaces to certain boundary
divisors). Thus, $G_x$ acts on $N_{G \cdot x/X,x}$ via a homomorphism 
$$
\rho_x: G_x \to \bG_m^r
$$ 
which is surjective by (iii). So $X$ is log homogeneous under $G$ by
Proposition \ref{prop:crit}.
\end{proof}

The converse of this corollary also fails in general. Indeed, all
$G$-orbit closures in a regular $G$-variety are also $G$-regular, and
hence nonsingular. But there exist log parallelizable varieties in
which most orbit closures are singular, see Example \ref{ex:tw}.

However, any {\it complete} log homogeneous variety $X$ under a
{\it reductive} group $G$ is regular by \cite[Prop.~2.5]{BB96}.
This will be generalized to certain non-complete varieties $X$ (and
reductive groups $G$) in Lemma \ref{lem:reg}, and then to arbitrary
groups $G$ (and complete varieties $X$) in Corollary \ref{cor:lhr}.

\subsection{Induced actions}
\label{subsec:ia}

Recall that a $G$-variety $X$ is {\it induced} from a homogeneous space
$G/H$ if there exists a $G$-equivariant morphism 
$f: X \to G/H$. Equivalently, 
$$
X \cong G \times^H Y,
$$ 
where $Y$ denotes the fiber of $f$ at the base point, and 
$G\times^H Y$ stands for the quotient of $G \times Y$ by the action of
$H$ via $h \cdot (g,y) = (gh^{-1},h \cdot y)$. We now show that
induction preserves homogeneity or parallelizability:

\begin{proposition}\label{prop:Gf}
Let $X$ be a $G$-variety and $D \subset X$ a reduced $G$-stable
divisor. Assume that $X = G \times^H Y$ for some closed subgroup 
$H \subseteq G$ and some closed $H$-stable subvariety  
$Y \subseteq X$, and let $E := D \cap Y$ (this is a reduced $H$-stable
divisor in $Y$). Then the following conditions are equivalent:

\smallskip

\noindent
{\rm (i)} $(X,D)$ is a homogeneous (resp.~parallelizable)
$G$-pair. 

\smallskip

\noindent
{\rm (ii)} $(Y,E)$ is a homogeneous (resp.~parallelizable)
$H^0$-pair. 
\end{proposition}

\begin{proof}
Clearly, $X$ is nonsingular if and only if so is $Y$. Consider then the
morphism 
$$
f: G \times Y \to X, \quad (g,y) \mapsto g \cdot y.
$$ 
Since $f$ is smooth, $D$ has normal crossings in $X$ if and only if
so has $f^{-1}(D)$ in $G \times Y$. As $f^{-1}(D) = G \times E$, the
latter condition means that $E$ has normal crossings in $Y$.

Next, consider the exact sequence of $H$-linearized sheaves on $Y$:
$$
0 \to \cT_Y \to \cT_X\vert_Y \to \cN_{Y/X} \to 0.
$$
The normal sheaf $\cN_{Y/X}$ is isomorphic as a $H$-linearized sheaf
to $\cO_Y \otimes \fg/\fh$, where $\fh$ denotes the Lie algebra of
$H$, and $\fg/\fh$ is regarded as a $H$-module. The sequence
$$
0 \to \cT_Y(- \log E) \to \cT_X(- \log D)\vert_Y \to \cN_{Y/X} \to 0 
$$
is also exact, since the composition 
$\cO_Y \otimes \fg \to \cT_X(- \log D)\vert_Y \to  \cO_Y \otimes \fg/\fh$
is surjective. This yields a commutative diagram of exact sequences
$$
\CD
0 @>>> \cO_Y \otimes \fh @>>> \cO_Y \otimes \fg @>>> 
\cO_Y \otimes \fg/\fh @>>> 0 \\
& & @VVV @VVV @VVV \\
0 @>>> \cT_Y( \log E) @>>> \cT_X(- \log D)\vert_Y @>>> \cN_{Y/X} @>>> 0, \\
\endCD
$$
where the right vertical map is an isomorphism. Thus, the surjectivity
(resp.~bijectivity) of the middle vertical map is equivalent to that
of the left vertical map.
\end{proof}

As a first application of this lemma, we construct a log
parallelizable variety which is not regular:

\begin{example}\label{ex:tw}
Consider the group $G := \GL_n$, where $n \ge 2$, and the 
subgroup $H$ of matrices having exactly one non-zero coefficient on
each line and each column; in other words, $H$ is the semidirect
product of the maximal torus $T$ of diagonal invertible matrices with
the symmetric group $S_n$ of permutation matrices, where $S_n$ acts on
$T$ by permuting the diagonal entries. Let $Y$ be the
affine $n$-space on which $H$ acts linearly via its standard
representation. Then $Y$ is log parallelizable under $T = H^0$, and
its boundary $E$ (the union of the coordinate hyperplanes
$E_1,\ldots,E_n$) is $H$-stable. The $H$-orbit closures in $E$ are
exactly the subsets $H \cdot (E_1 \cap \cdots \cap E_r)$, where 
$r = 1,\ldots,n$. In particular, $E$ is the closure of a unique
$H$-orbit.

Since the $H$-action on $Y$ extends to a $G$-action, 
the induced variety $X := G \times^H Y$ is just the product 
$G/H \times Y$ on which $G$ acts diagonally; since $H$ is reductive,
$X$ is affine. By the above lemma, $X$ is log parallelizable under
$G$; its boundary  $D := G \times^H E$ is a $G$-orbit closure, and
hence is irreducible. More generally, the $G$-orbit closures in
$X$ are the subsets
$$
X^r := G \cdot (E_1 \cap \cdots \cap E_r),
$$ 
where $r=1, \ldots, n$, together with $X$. 

We claim that the singular locus of each $X^r$ is $X^{r+1}$. This
claim is a consequence of Corollary \ref{cor:fo}, or may be seen
directly as follows. Let $H_r$ be the stabilizer in $H$ of the subset 
$E_1 \cap \cdots \cap E_r$, so that $H_r = T S_{n-r}$. Then the
natural map
$$
H \times^{H_r} (E_1 \cap \cdots \cap E_r) \to 
H \cdot (E_1 \cap \cdots \cap E_r)
$$
is finite and restricts to an isomorphism over the complement of 
$H \cdot (E_1 \cap \cdots \cap E_{r+1})$, but over no larger open
subset. It follows that the natural map
$$
G \times^{H_r} (E_1 \cap \cdots \cap E_r) \to X^r
$$
is finite and birational with exceptional set $X^{r+1}$. Since 
$X^r \setminus X^{r+1}$ is a unique orbit, this implies our claim.

As a consequence, $X$ is not $G$-regular. This may also be seen by
considering the base point $x := H \cdot (1,0)$ of the closed
orbit. Then $G_x = H$ acts on the normal space 
$N_{G \cdot x/X,x} \cong \bA^n$ via its standard representation; in
particular, $G_x$ permutes transitively the coordinate lines
$L_1, \ldots, L_n$.
\end{example}

As another application of Proposition \ref{prop:Gf}, we show that 
homogeneity and parallelizability are also preserved by taking general 
fibers: 

\begin{corollary}\label{cor:lhf}
Consider a homogeneous (resp.~parallelizable) $G$-pair $(X,D)$, a
variety $X'$, and a proper morphism $f : X \to X'$. Let $Y$  
be a connected component of the fiber of $f$ at a point of the open
orbit $X_0 = X \setminus D$. Then $(Y, D \cap Y)$ is a homogeneous
(resp.~parallelizable) pair.  
\end{corollary} 

\begin{proof}
Using the Stein factorization, we may assume that the natural map
$\cO_{X'} \to f_* \cO_X$ is an isomorphism; then all the fibers of $f$
are connected, and the general fibers are irreducible. 

By a result of Blanchard (see \cite[Sec.~2.4]{Ak95} for a proof in the
setting of complex geometry, which may be adapted readily to our
algebraic setting), the $G$-action on $X$ descends to an action on
$X'$ such that $f$ is equivariant. Let $X'_0 \cong G/H'$ be the open
$G$-orbit in $X'$. Then $f^{-1}(X'_0)= f^{-1}f(X_0)$ is an
open $G$-stable subset of $X$, equivariantly isomorphic to 
$G \times^{H'} Y$. In particular, all fibers over $X'_0$ are
isomorphic to $Y$, and hence the latter is a nonsingular variety. So
the assertions follow from Proposition \ref{prop:Gf}.
\end{proof}

\subsection{Local models}
\label{subsec:lm}

A simple example of a log parallelizable variety is the affine space
$\bA^r$ with its standard action of $\bG_m^r$; one may also construct
induced versions $G \times^H \bA^r$, where $H$ is a subtorus of $G$,
acting on $\bA^r$ via an isogeny to $\bG_m^r$. We now show that this
construction yields local models for all log homogeneous varieties: 

\begin{proposition}\label{prop:lm}
Let $(X,D)$ be a homogeneous $G$-pair. Choose $x \in X$ and let 
$\rho_x: G_x^0 \to \bG_m^r$ be as in (\ref{eqn:rho}). Then 
there exist a subtorus $H \subseteq G_x^0$ and a $H$-stable locally
closed subvariety $Y \subseteq X$ such that:

\smallskip

\noindent
{\rm (i)} The restriction of $\rho_x$ to $H$ is an isogeny to
$\bG_m^r$. 

\smallskip

\noindent
{\rm (ii)} $Y \cong \bA^r$ on which $H$ acts via the restriction 
$\rho_x : H \to \bG_m^r$ and the standard action of $\bG_m^r$ on
$\bA^r$.

\smallskip

\noindent
{\rm (iii)} The natural map $p: G \times^H Y \to X$ is smooth and
$p^{-1}(D) \cong G \times^H E$, where $E \subset \bA^r$
denotes the union of all coordinate hyperplanes.

\smallskip

Moreover, $(X,D)$ is parallelizable if and only if $p$ is \'etale.
\end{proposition}

\begin{proof}
We use some of the ingredients of the Luna slice theorem \cite{Lu73}.
Since $\rho_x$ is surjective, we may choose a subtorus 
$H \subseteq G_x^0$ satisfying (i). Then there exists a $H$-stable
decomposition 
$$
T_x X = T_x (G \cdot x) \oplus N = \fg \cdot x \oplus N,
$$
where the $H$-module $N$ is isomorphic to $N_{G\cdot x/X,x}$.
Moreover, the $H$-fixed point $x$ admits an affine $H$-stable
neighborhood $X_x \subseteq X$ together with an $H$-equivariant map
$$
\varphi: X_x \to T_x X
$$
such that the differential of $\varphi$ at $x$ is the identity map of
$T_x X$. In particular, $\varphi$ is \'etale at $x$.

Let $Y := \varphi^{-1}(N)$. This is an affine $H$-stable locally
closed subvariety of $X$, containing $x$ and nonsingular at that point.
Moreover, $T_x Y \cong N_{G\cdot x/X,x}$ as $H$-modules. By the graded
Nakayama lemma, it follows that $Y$ satisfies (ii). For (iii), note
that $p$ is smooth at the point 
$H \cdot (1,x)$ of $G\times^H Y$, since $Y$ is a slice to the orbit
$G\cdot x$ at $x$. Thus, $p$ is smooth in a $G$-stable
neighborhood of $H \cdot (1,x)$, which must be the whole 
$G\times^H Y$.

Finally, $p$ is \'etale if and only if $H = G_x^0$; equivalently,
$\rho_x$ is an isogeny. By Proposition \ref{prop:crit}, this means
that $(X,D)$ is parallelizable.
\end{proof}

From the above proposition, we deduce that homogeneity and
parallelizability are preserved under equivariant blowing up:

\begin{proposition}\label{prop:mor}
{\rm (i)} Let $(X,D)$ be a homogeneous (resp.~parallelizable)
$G$-pair, $X'$ a nonsingular $G$-variety, $f:X' \to X$ a birational
$G$-equivariant morphism, and $D'$ the reduced inverse image of
$D$. Then $(X',D')$ is a homogeneous (resp.~parallelizable) $G$-pair.
\smallskip

\noindent
{\rm (ii)} As a partial converse, given two $G$-pairs $(X,D)$, 
$(X',D')$ and a surjective birational $G$-equivariant morphism 
$f : X' \to X$ such that $D'$ is the reduced inverse image of
$D$, if $(X',D')$ is $G$-parallelizable, then so is $(X,D)$.
\end{proposition}

\begin{proof}
(i) Choose $x \in X$ and let $H$, $Y$ be as in Proposition
\ref{prop:lm}. Form the cartesian square
$$
\CD
W  @>{\varphi}>> G \times^H Y \\
@V{\pi}VV        @V{p}VV \\
X' @>{f}>>       X \\
\endCD
$$
Then $W = G \times^H Y'$, where $Y' := f^{-1}(Y)$. Since $p$ is smooth,
then so is $\pi$. Thus, $Y'$ is a nonsingular variety; it is also
toric under $\bG_m^r$. Hence $Y'$ is log parallelizable by
\cite[Prop.~3.1]{Od88}; its boundary $\partial Y'$ equals the reduced
inverse image $f^{-1}(D \cap Y)_{\red}$. It follows that the
$G$-variety $W$ is log parallelizable with boundary
$$
\partial W = G \times^H \partial Y'= \varphi^{-1} p^{-1}(D)_{\red}. 
$$
Since $\pi$ is smooth and $\partial W = \pi^{-1}(D')$, the divisor
$D'$ has normal crossings. Moreover, the surjective map 
$d \pi: \cT_W \to \pi^* \cT_{X'}$ 
restricts to a surjective map
$$
\cT_W( - \log \partial W) \to \pi^* \cT_{X'}(- \log D').
$$
It follows that $(X',D')$ is $G$-homogeneous (resp.~parallelizable if
so is $(X,D)$).

(ii) The assumption implies that the map $op_{X,D}$ of (\ref{eqn:opd})
is an isomorphism in codimension $1$. Since this is a map between
locally free sheaves of the same rank, it is an isomorphism.
\end{proof}

\begin{remark}
The above statement (ii) does not extend to homogeneous pairs. For
example, let $G := \PGL_{n+1}$ acting diagonally on 
$X := \bP^n \times \bI_n$, where $\bI_n$ denotes the incidence variety
consisting of pairs $(p,h)$ such that $p \in \bP^n$ and 
$h \subset \bP^n$ is a hyperplane through $p$. Then the nonsingular
$G$-variety $X$ contains an open $G$-orbit with complement $D$
consisting of those triples $(p_1,p_2,h)$ such that 
$p_1,p_2 \in h$. Thus, $D$ is a nonsingular prime divisor in $X$,
consisting of two $G$-orbits; the closed orbit $Y$ (where $p_1 = p_2$)
is isomorphic to $\bI_n$. Hence $X$ consists of two strata but three
orbits.

In particular, $(X,D)$ is a nonhomogeneous $G$-pair. But one may check
that the blowing up of $Y$ in $X$ is log homogeneous under $G$.
\end{remark}

\subsection{The Albanese fibration}
\label{subsec:af2}

For any compact K\"ahler almost homogenous manifold, one knows that
the Albanese map is an equivariant fibration with connected fibers;
see e.g.~\cite[Sec.~3.9]{Ak95}. The following result is an algebraic
analogue for varieties that need not be complete:

\begin{proposition}\label{prop:ar}
Let $X$ be a nonsingular almost homogeneous $G$-variety with open
orbit $X_0$ and let $\cA(X_0) = G/I(X_0)$ be the Albanese variety
(as defined in Subsec.~\ref{subsec:af1}). Then:

\smallskip

\noindent
{\rm (i)} The Albanese morphism
$$
\alpha : X  \to \cA(X)
$$
is a $G$-equivariant fibration, and $\cA(X) = \cA(X_0)$. In particular, 
$$
X \cong G \times^I Y,
$$ 
where $I := I(X_0)$ is a closed subgroup of $G$ containing $G_{\aff}$,
and the fiber $Y$ of $\alpha$ is a nonsingular $I$-variety, almost
homogeneous under $G_{\aff}$.

\smallskip

\noindent
{\rm (ii)} If $G$ acts faithfully on $X$, then $I$ (and hence
$G_{\aff}$) acts faithfully on $Y$.

\smallskip

\noindent
{\rm (iii)} Given a reduced $G$-stable divisor $D$ in $X$, the pair
$(X,D)$ is homogeneous (resp.~parallelizable) under $G$ if and only if 
$(Y,D \cap Y)$ is homogeneous (resp.~parallelizable) under
$G_{\aff}$. Then the restriction to $\cT_X(- \log D)$ of the natural
map 
$$
\cT_X \to \alpha^* \cT_{\cA(X)} \cong \cO_X \otimes \fg/\fg_{\aff}
$$
fits into an exact sequence 
$$
0 \to \cT_{\alpha}( - \log D) \to \cT_X(- \log D) \to 
\cO_X \otimes \fg/\fg_{\aff} \to 0,
$$
where $\cT_{\alpha}( - \log D)$ is a $G$-linearized locally free sheaf,
generated by the image of $\fg_{\aff}$. Moreover,
$\cT_{\alpha}( - \log D)\vert_Y \cong \cT_Y( - \log (D \cap Y))$.
\end{proposition}

\begin{proof}
(i) By \cite[Thm.~3.1]{Mi86}, the Albanese map $X_0 \to \cA(X_0)$
extends to a morphism $\alpha: X \to \cA(X_0)$. Clearly, $\alpha$ is a 
$G$-equivariant fibration with fibers being almost homogeneous
$G_{\aff}$-varieties. In particular, any morphism from a fiber to an
abelian variety is constant; it follows that $\alpha$ is the Albanese 
morphism of $X$. 

By (\ref{eqn:gz}), $X = C^0 \cdot Y$ which implies (ii).

(iii) The equivalences follow from (i) together with Proposition
\ref{prop:Gf}. The final assertions are consequences of the
commutative diagram
$$
\CD
0 @>>> \cO_X \otimes \fg_{\aff} @>>> \cO_X \otimes \fg @>>> 
\cO_X \otimes \fg/\fg_{\aff} @>>> 0 \\
& & @VVV @VVV @V{\cong}VV \\
0 @>>> \cT_{\alpha}( \log D) @>>> \cT_X(- \log D) @>>> 
\alpha^*\cT_{\cA(X)} @>>> 0. \\
\endCD
$$
\end{proof}

Using this Albanese fibration, we reduce the description of the
$G$-equivariant automorphism group $\Aut_G(X)$ to that of $\Aut_I(Y)$: 

\begin{lemma}\label{lem:ea}
With the notation and assumptions of the above proposition, 
the group $\Aut_G(X)$ is algebraic and fits into an exact sequence of
such groups
\begin{equation}\label{eqn:ig}
1 \to \Aut_I(Y) \to \Aut_G(X) \to \cA(X) \to 0.
\end{equation}
Moreover, $\Aut_I(X)$ is affine.
\end{lemma}

\begin{proof}
We may identify $\Aut_G(X)$ with a subgroup of 
$\Aut_G(X_0) \cong N_G(H)/H$, where $X_0 \cong G/H$. 
To show that $\Aut_G(X)$ is algebraic, we check that
it is closed in $N_G(H)/H$. For this, consider the ``graph''
$$
Y_0 := \{(x,g \cdot x,g)~\vert~ x \in X_0, ~g \in N_G(H)/H \}
\subseteq X_0 \times X_0 \times N_G(H)/H
$$
and its closure $Y$ in $X \times X \times N_G(H)/H$. Then the
projection 
$$
p_1 \times p_3: Y \to X \times N_G(H)/H
$$ 
is an isomorphism over the dense open subset $X_0 \times N_G(H)/H$.
Let $E \subset Y$ be the exceptional set of $p_1 \times p_3$. Then
$\Aut_G(X)$ consists of those $\gamma \in N_G(H)/H$ such that
$X \times X \times \{ \gamma \}$ does not meet $E$. Thus, $\Aut_G(X)$
is a constructible subgroup of $N_G(H)/H$, and hence is closed.

Next, we obtain the exact sequence (\ref{eqn:ig}). By the universal
property of the Albanese morphism, any $G$-equivariant automorphism of
$X$ induces an equivariant automorphism of the abelian variety
$\cA(X)$, that is, a translation. This yields a group homomorphism
$$
f: \Aut_G(X) \to \cA(X)
$$
which is the restriction of the natural homomorphism
$$
\Aut_G(X_0)= N_G(H)/H \to G/G_{\aff}H = \cA(G/H) = \cA(X)
$$
and hence is algebraic. The composition 
$C^0 \to \Aut_G(X) \to \cA(X)$ is surjective by Lemma \ref{lem:G},
so that $f$ is surjective. Since $X \cong G \times^I Y$, the kernel of 
$f$ is isomorphic to $\Aut_I(Y)$.

As in the first step of the proof, $\Aut_I(Y)$ is a closed
subgroup of $N_I(H)/H$, and hence is affine since $I$ is.
\end{proof}

\begin{remarks}\label{rem:sim}
We consider again a nonsingular variety $X$, almost homogeneous under 
a faithful action on $G$, and we use the notation of Proposition
\ref{prop:ar}. 

\smallskip

\noindent
(i) Let $I = G_{\aff} F$ be as in (\ref{eqn:if}). Then $F$ acts on $Y$
via a homomorphism 
$$
\varphi : F \to \Aut_I(Y) \subseteq \Aut_{G_{\aff}}(Y).
$$
Thus, if every $G_{\aff}$-equivariant automorphism of $Y$ comes from
some element of $G_{\aff}$, then $I = G_{\aff}$ (since $I$ acts
faithfully on $Y$). Equivalently, the isogeny $\cA(G) \to \cA(X)$ is
an isomorphism under this assumption.

\smallskip

\noindent
(ii) In general, we may reduce to the case where $\cA(G) = \cA(X)$,
just by replacing $X$ with the finite \'etale cover
$G \times^{G_{\aff}} Y$ and keeping $G$ unchanged.

We may reduce further to the case where the isogeny 
$\cA(C^0) \to \cA(G)$ (Lemma \ref{lem:G}) is an isomorphism, by
replacing $X$ with the finite \'etale cover 
$$
X' := C^0 \times^{C^0_{\aff}} Y
$$ 
and $G$ with the finite cover
$$
G' := C^0 \times^{C^0_{\aff}} G_{\aff}.
$$ 
Then the Albanese map
$$
\alpha' : C^0 \times^{C^0_{\aff}} Y \to C^0/C^0_{\aff}
$$
is locally trivial for the Zariski topology, since the algebraic group
$C^0_{\aff}$ is commutative, connected and affine.
\end{remarks}

In the above remark, both covers of $X$ are just the
pull-backs under the Albanese map of the isogenies 
$\cA(G) \to \cA(X)$, resp.~$\cA(C^0) \to \cA(G)$. We now show that
every finite \'etale cover of an almost homogeneous variety is
obtained in that way, under the assumption of completeness:

\begin{proposition}\label{prop:ec}
Let $X$ be a complete nonsingular variety, almost homogeneous under a
faithful $G$-action. Let $f: X' \to X$ be an \'etale morphism, where
$X'$ is a complete (nonsingular) variety. Then $X'$ is almost
homogeneous under a faithful action of a finite cover $G'$ of
$G$. Moreover, $f$ induces an isogeny $\varphi: \cA(X') \to \cA(X)$,
and the square 
\begin{equation}\label{eqn:sq}
\CD
X' @>{\alpha'}>> \cA(X') \\
@V{f}VV @V{\varphi}VV \\
X @>{\alpha}>> \cA(X)
\endCD
\end{equation}
is cartesian.
\end{proposition}

\begin{proof}
Note that any derivation of $\cO_X$ extends to a unique derivation
of $f_*\cO_{X'}$. This yields an injective homomorphism of Lie algebras
$$
i : \Gamma(X,\cT_X) \to \Gamma(X',\cT_{X'}).
$$ 
Moreover, the image of $i$ is an algebraic Lie subalgebra, since it
consists of those derivations that preserve the field of rational
functions on $X$. Thus, $i(\fg)$ is algebraic as well; let 
$G' \subseteq \Aut^0(X')$ be the associated algebraic group. Then $G'$
is a finite cover of $G$, and has an open orbit in $X'$.

The universal property of the Albanese morphism yields a
commutative square (\ref{eqn:sq}), where $\varphi$ is a
homomorphism. We also have a commutative square
$$
\CD
\cA(G') @>>> \cA(X')\\ 
@VVV @V{\varphi}VV \\
\cA(G) @>>> \cA(X),
\endCD
$$ 
where the three nonlabeled arrows are isogenies. It follows that
$\varphi$ is an isogeny as well. So the fibers $Y$, $Y'$ of $\alpha$,
resp.~$\alpha'$, have the same dimension and satisfy 
$Y' \subseteq f^{-1}(Y)$. But $Y$ is unirational and hence
simply connected (see \cite[Exp.~XI, Cor.~1.3]{SGA1}); thus, $f$
restricts to an isomorphism $Y' \to Y$. This implies that the square
(\ref{eqn:sq}) is cartesian.
\end{proof}

\subsection{The complete log parallelizable case}
\label{subsec:lp}

We shall obtain an algebraic version of the main result of \cite{Wi04} 
describing the log parallelizable compact ``weakly K\"ahler''
manifolds. To state our version, recall that a
{\it semiabelian variety} is an algebraic group obtained as an
extension of an abelian variety by a torus; any semiabelian variety is 
connected and commutative.

\begin{theorem}\label{thm:lpt}
Let $X$ be a complete nonsingular variety on which a connected
algebraic group $G$ acts faithfully. Then the following conditions are
equivalent:

\smallskip

\noindent
{\rm (i)} $X$ is log parallelizable under $G$.

\smallskip

\noindent
{\rm (ii)} $X \cong G \times^{G_{\aff}} Y$, where $G_{\aff}$ is a torus
and $Y$ is a complete nonsingular toric variety under $G_{\aff}$.

\smallskip

\noindent
{\rm (iii)} $G$ is a semiabelian variety and has an open orbit in $X$.

\smallskip

Under one of these conditions, $G = \Aut^0(X,D)$, where $D$ denotes
the boundary.
\end{theorem}

\begin{proof}
(i)$\Rightarrow$(ii) We use the notation of Proposition \ref{prop:ar}.
Note that $Y$ is log homogeneous under $G_{\aff}$. Choose $y \in Y$
such that the orbit $G_{\aff} \cdot y$ is closed. Then the isotropy
group $G_{\aff,y}$ is a parabolic subgroup of $G_{\aff}$. Hence
$G_{\aff,y}$ is connected and contains a maximal unipotent subgroup of
$G_{\aff}$. On the other hand, $G_y^0$ is contained in $G_{\aff}$, and
is isogenous to a torus by Proposition \ref{prop:crit}. It follows
that $G_{\aff,y}$ is a torus and, in turn, that $G_{\aff}$ is a torus
as well. 

Hence $Y$ is a toric $G_{\aff}$-variety, so that the natural map 
$G_{\aff} \to \Aut_{G_{\aff}}(Y)$ is an isomorphism. This implies
that $I = G_{\aff}$ by Remark \ref{rem:sim}(i).

(ii)$\Rightarrow$(iii) is obvious.

(iii)$\Rightarrow$(i) Since $G$ is commutative and acts faithfully in
$X$, its open orbit is isomorphic to $G$ itself. Together with
Proposition \ref{prop:ar}, it follows that $X = G \times^{G_{\aff}} Y$,
where $Y$ contains an open orbit of the torus $G_{\aff}$. Thus, $Y$ is
log parallelizable (see e.g. \cite[Prop.~3.1]{Od88}) and hence so is 
$X$.

Under the condition (i), we have the equalities of Lie algebras
$$
\fg = \Gamma(X,\cO_X \otimes \fg) = \Gamma(X,\cT_X(-\log D))
$$
and hence $G = \Aut^0(X,D)$. 
\end{proof}

This displays the close relationship between log parallelizable and
{\it semiabelic} varieties; recall from \cite{Al02} that the latter
are the normal varieties on which a semiabelian variety acts with
finitely many orbits, such that all isotropy groups are tori. By the
above theorem, being log parallelizable or semiabelic are equivalent
for complete nonsingular varieties.

\section{Complete log homogeneous varieties}
\label{sec:complete}

\subsection{Relation to spherical varieties}
\label{subsec:rsv}

We consider a complete log homogeneous $G$-variety $X$ and write 
$X \cong G\times^I Y$ as in Proposition \ref{prop:ar}, so that $Y$ is
a complete nonsingular variety, log homogeneous under $G_{\aff}$. In 
this subsection, we obtain several preliminary results about the
$G_{\aff}$-variety $Y$; together, they will imply our main structure
theorem in the next subsection.

Choose a Levi subgroup $L \subseteq G_{\aff}$, i.e., a maximal closed 
connected reductive subgroup. Recall the semidirect product
decomposition $G_{\aff} = R_u(G_{\aff}) L$, where $R_u(G_{\aff})$
denotes the unipotent radical; moreover, any two Levi subgroups are
conjugate by an element of $R_u(G_{\aff})$. 

\begin{lemma}\label{lem:lhp}
For a $G_{\aff}$-pair $(Y,E)$ where $Y$ is complete, the following
conditions are equivalent: 

\smallskip

\noindent
{\rm (i)} $(Y,E)$ is homogeneous under $G_{\aff}$.

\smallskip

\noindent
{\rm (ii)} $(Y,E)$ is homogeneous under $L$ along each closed
$G_{\aff}$-orbit. 

\smallskip

If one of these conditions holds, then there exists a smallest open
$L$-stable subset $Y_L \subseteq Y$ which contains every closed
$G_{\aff}$-orbit. Then the pair $(Y_L, E \cap Y_L)$ is homogeneous
under $L$, and every $G_{\aff}$-orbit in $Y$ meets $Y_L$ along a
unique $L$-orbit. In particular, the closed $L$-orbits in $Y_L$ are
exactly the closed $G_{\aff}$-orbits in $Y$.
\end{lemma}

\begin{proof}
(i)$\Rightarrow$(ii) As in the proof of Theorem \ref{thm:lpt},
choose $y \in Y$ such that the orbit $G \cdot y$ is closed, i.e.,
$G_y$ is a parabolic subgroup of $G_{\aff}$. Then 
$G_{\aff,y} = R_u(G_{\aff})  L_y = G_{\aff,y}^0$ 
and $G_{\aff} \cdot y = L \cdot y$. Clearly, the homomorphism 
$\rho_y: G_y^0 \to \bG_m^r$ of (\ref{eqn:rho}) has a trivial
restriction to $R_u(G_{\aff})$. Thus, by Proposition \ref{prop:crit},
$(Y,E)$ is homogeneous under $L$ along $G_{\aff} \cdot y$. 

(ii)$\Rightarrow$(i) The assumption implies that the map $op_{Y,E}$ is
surjective along any closed $G_{\aff}$-orbit. Thus, $op_{Y,E}$ is
surjective everywhere.

If (ii) holds, then there exists an open $L$-stable subvariety
$Y' \subseteq Y$ containing all the closed $G_{\aff}$-orbits, such
that the pair $(Y',E \cap Y')$ is homogeneous under $L$. In
particular, $Y'$ contains only finitely many $L$-orbits, and every
closed $G_{\aff}$-orbit in $Y$ is a closed $L$-orbit in $Y'$. Let
$Y_L$ be the set of those $y \in Y$ such that the orbit closure
$\overline{L \cdot y}$ contains a closed $G_{\aff}$-orbit. Then $Y_L$
is contained in $Y'$ and is the smallest common open $L$-stable
neighborhood of all the closed $G_{\aff}$-orbits. 

As a consequence $G_{\aff} \cdot Y_L$ is an open $G$-stable
neighborhood of the closed $G$-orbits, so that 
$G_{\aff} \cdot Y_L = Y$. In other words, 
every $G_{\aff}$-orbit $\Omega$ in $Y$ meets $Y_L$. Since the
$L$-orbits in $Y_L$ are the strata of the pair $(Y_L, E \cap Y_L)$,
i.e., the intersections of $Y_L$ with the strata of $(Y,E)$, it
follows that $\Omega \cap Y_L$ is a unique $L$-orbit.
\end{proof}

Next, we obtain a partial converse of Corollary \ref{cor:rlh}:

\begin{lemma}\label{lem:reg}
Let $Y$ be a log homogeneous variety under a connected reductive group
$L$. If every closed orbit is complete, then the $L$-variety $Y$ is
spherical and regular. Moreover, its closed orbits are all isomorphic.
\end{lemma}

\begin{proof}
This is proved in \cite[Prop.~2.2.1]{BB96} under the assumption that
$Y$ is complete; the argument there may be adapted as follows.
Choose $y\in Y$ such that the orbit $L \cdot y$ is closed, so that
$L_y$ is a parabolic subgroup of $L$. Let $P \subseteq L$ be a
parabolic subgroup opposite to $L_y$, so that
$$
M := P \cap L_y
$$ 
is a Levi subgroup of both $P$ and $L_y$. By the local structure
theorem (see e.g. \cite[Sec.~2.1]{Kn94a}), there exists a locally
closed affine $M$-stable subvariety $S \subseteq Y$ containing $y$,
such that the map 
$$
R_u(P) \times S \to Y, \quad (g,s) \mapsto g \cdot s
$$ 
is an open immersion. In particular, $S$ is nonsingular and meets
transversally $L \cdot y$ at the unique point $y$. This yields a
$M$-equivariant isomorphism 
$$
T_y S \cong N_{L\cdot y/Y,y}.
$$
In particular, $M$ acts on $T_y S \cong \bA^r$ via a surjective
homomorphism to $\bG_m^r$, and hence the derived subgroup of $M$ fixes
$T_y S$ pointwise. As in the proof of Proposition \ref{prop:lm}, it
follows that $S \cong \bA^r$ as $M$-varieties, and 
the boundary of $Y$ meets $S$ along the union of all the coordinate
hyperplanes in $\bA^r$. Thus, the $L$-variety $Y$ is regular in a
neighborhood of $L \cdot y$. Moreover, the image of 
$R_u(P) \times \bG_m^r$ in $Y$ is the open orbit of a Borel subgroup
$B$ of $L$, so that $Y$ is spherical as well. Finally, one checks that
the stabilizer in $L$ of this open $B$-orbit equals the parabolic
subgroup $P$. In particular, the conjugacy class of $P$ is independent
of the closed orbit $L \cdot y$. Thus, the same holds for the
conjugacy class of the opposite parabolic subgroup $L_y$.
\end{proof}

Finally, we record the following result, a consequence of
\cite[Thm.~5.1, Cor.~5.6]{Kn96}.

\begin{lemma}\label{lem:aut}
Let $Y$ be a spherical $L$-variety, and $Y' \subset Y$ an open
$L$-stable subset. Then $\Aut_L(Y)$ is a diagonalizable algebraic
group and preserves $Y'$; in particular, $\Aut^0_L(Y)$ is a torus. 
Moreover, the restriction map 
$\Aut^0_L(Y) \to \Aut^0_L(Y')$ is an isomorphism.
\end{lemma}

\subsection{Structure}
\label{subsec:str}

We now come to our main result:

\begin{theorem}\label{thm:str}
Let $G$ be a connected algebraic group and choose a Levi subgroup
$L\subseteq G_{\aff}$. Then any complete log homogeneous $G$-variety
$X$ may be written uniquely as $G \times^I Y$, where

\smallskip

\noindent
{\rm (i)} $I \subseteq G$ is a closed subgroup containing $G_{\aff}$ as
a subgroup of finite index.

\smallskip

\noindent
{\rm (ii)} $Y$ is a complete nonsingular $I$-variety containing an
open $L$-stable subset $Y_L$ such that 

\smallskip

{\rm (a)} the $L$-variety $Y_L$ is regular, and  

\smallskip

{\rm (b)} every $G_{\aff}$-orbit in $Y$ meets $Y_L$ along a unique
$L$-orbit.

\smallskip

In particular, the $L$-variety $Y$ is spherical, and the projection 
$X \to G/I$ is the Albanese morphism. 

Conversely, given $I$ and $Y$ satisfying {\rm (i)} and {\rm (ii)}, the 
$G$-variety $X := G \times^I Y$ is log homogeneous. Moreover, each
$I$-orbit in $Y$ is a unique $G_{\aff}$-orbit; in particular, each
$G$-orbit in $X$ meets $Y$ (resp.~$Y_L$) along a unique orbit of
$G_{\aff}$ (resp.~$L$).
\end{theorem}

\begin{proof}
All the assertions follow from Proposition \ref{prop:ar} and Lemmas
\ref{lem:lhp}, \ref{lem:reg}, except for the equality of $I$-orbits
and $G_{\aff}$-orbits. For the latter, let again $I = G_{\aff} F$ be
as in (\ref{eqn:if}). The finite subgroup $F \subseteq C^0$ 
acts on $Y$ by $G_{\aff}$-equivariant automorphisms. Thus, $F$
preserves all $L$-orbits by Lemma \ref{lem:aut}, and hence all
$G_{\aff}$-orbits. 
\end{proof}

Next, we obtain some remarkable consequences of this structure theorem.

\begin{corollary}\label{cor:lhr}
Any complete log homogeneous $G$-variety $X$ is $G$-regular; in
particular, any $G$-stable subvariety is log homogeneous. Moreover,
the irreducible components of fibers of any morphism $f : X \to X'$
are log homogeneous.
\end{corollary}

\begin{proof}
The first assertion follows from Theorem \ref{thm:str} together with
the $L$-regularity of $Y_L$ (Lemma \ref{lem:reg}). 

For the second assertion, let $C$ be an irreducible component of a
fiber and let $X'$ be the closure of $G \cdot C$ in $X$. Then $X'$ is
a $G$-stable subvariety of the regular variety $X$, and hence is
regular as well. So the assertion follows by applying Corollary
\ref{cor:lhf} to $f \vert_{X'}$.
\end{proof}

Also, we rephrase the most delicate condition (b) in the statement of
Theorem \ref{thm:str}, to make it easier to check:

\begin{lemma}\label{lem:char}
Let $Y$ be a complete nonsingular $G_{\aff}$-variety, and
$Y_L \subseteq Y$ an $L$-stable open subset. Assume that $Y_L$ is
$L$-regular and contains all the closed $G_{\aff}$-orbits in $Y$. 
Let $\partial Y_L$ be the boundary of the $L$-variety $Y_L$, and 
$E := G_{\aff} \cdot \partial Y_L$. Then the following conditions are 
equivalent:

\noindent
{\rm (b)} Every $G_{\aff}$-orbit in $Y$ meets $Y_L$ along a unique
$L$-orbit.

\noindent
{\rm (c)} $Y_L$ is not contained in $E$.

\noindent
{\rm (d)} $\partial Y_L$ is stable under $\fg_{\aff}$ (acting on $Y_L$ by
vector fields).

Under one of these conditions, $(Y,E)$ is a $G_{\aff}$-homogeneous pair. 
\end{lemma}

\begin{proof}
Notice that $G_{\aff} \cdot Y_L$ is a $G$-stable open subset of $Y$
containing all the closed $G$-orbits, whence $G_{\aff} \cdot Y_L = Y$.
In other words, every $G_{\aff}$-orbit in $Y$ meets $Y_L$.

(b)$\Rightarrow$(c) is obvious.

(c)$\Rightarrow$(d) The assumption implies the equality
$Y_L \cap E = \partial Y_L$, which in turn
implies the desired statement. 

(d)$\Rightarrow$(b) By our assumption, $E$ does not meet the open
$G_{\aff}$-orbit in $Y$. It follows that $E \cap Y_L = \partial Y_L$. 
On the other hand, $E$ is a divisor with normal crossings (since so is
$\partial Y_L$), and its complement is the open $G_{\aff}$-orbit in
$Y$. Moreover, each stratum of $(Y,E)$ meets $Y_L$ along a unique
stratum of $(Y_L,\partial Y_L)$, i.e., a unique $L$-orbit. Together
with Proposition \ref{prop:crit}, this completes the proof.
\end{proof}

\subsection{The Tits morphism}
\label{subsec:tm}

In this subsection, we consider a complete log homogeneous variety
$X$ under a faithful action of $G$. Let 
$$
\tau : X \to \cL, \quad x \mapsto \fg_{(x)}
$$
be the $G$-morphism defined in Corollary \ref{cor:fo}. We say that
$\tau$ is the {\it Tits morphism} of $X$, and we denote its image by 
$\cL(X) \subseteq \cL$. 

Clearly, $\cL(X)$ is a projective $G$-variety, fixed pointwise by $C$.
Moreover, in the notation of Theorem \ref{thm:str}, $\cL(X) = \cL(Y)$,
and the restriction $\tau \vert_Y$ is the Tits morphism of the
complete log homogeneous $G_{\aff}$-variety $Y$. The sheaf 
$\cT_X ( - \log D)$ is the pull-back under $\tau$ of the quotient
sheaf on $\cL(X)$ (regarded as a subvariety of the Grassmannian of
subspaces of $\fg$). Further properties of $\tau$ are gathered in the
following:    

\begin{proposition}\label{prop:tm}
{\rm (i)} The $G$-variety $\cL(X)$ contains a unique closed orbit
$\cF(X)$, and $\tau$ maps every closed $G_{\aff}$-orbit in $Y$
isomorphically to $\cF(X)$.  

\smallskip

\noindent
{\rm (ii)} The irreducible components of the fibers of $\tau$ are
log parallelizable varieties.
\end{proposition}

\begin{proof}
(i) By the above observations, we may assume that $G$ is affine,
i.e., $X = Y$. Let $x \in X$ such that the orbit $G \cdot x$ is
closed. Then $G \cdot \tau(x) \cong G/N_G(\fg_{(x)})$. Moreover, the
Lie subalgebra $\fg_{(x)} \subseteq \fg_x$ is the intersection of the
kernels of differentials of characters of the isotropy group $G_x$, a
parabolic subgroup of $G$. It follows that $N_G(\fg_{(x)}) = G_x$;
hence $\tau$ restricts to an isomorphism
$$
G \cdot x \cong G \cdot \tau(x) = G \cdot \fg_x.
$$
Moreover, $G_x = R_u(G) L_x$, where $L$ is a Levi subgroup of
$G$ and the conjugacy class of $L_x$ in $L$ is 
independent of the closed $G$-orbit $G \cdot x$ by Lemma
\ref{lem:reg}. Thus, the orbit $G \cdot \fg_x$ is also independent of
the closed $G$-orbit.

(ii) By Corollary \ref{cor:lhr}, the irreducible components of fibers
of $\tau$ are nonsingular. To show that they are log parallelizable,
it suffices to consider fibers at points $y \in Y$. Then
$$
\tau^{-1}\tau(y) = C^0 (\tau^{-1}\tau(y) \cap Y)
$$
is mapped by $\alpha$ onto $\cA(X)= C^0/(C^0 \cap I)$. Thus,
\begin{equation}\label{eqn:fi}
\tau^{-1}\tau(y) \cong 
C^0 \times^{C^0 \cap I} (\tau^{-1}\tau(y) \cap Y).
\end{equation}
Since $\tau^{-1}\tau(y) \cap Y$ is a fiber of the Tits morphism
of $Y$, we are reduced to checking that the irreducible
components of fibers of that morphism are toric varieties. In other
words, we may assume again that $G$ is affine.

Choose a point $x \in \cF(X)$. Applying the local structure theorem as
in the proof of Lemma \ref{lem:reg}, we obtain a parabolic subgroup
$P \subseteq L$ and a locally closed affine
subvariety $S \subseteq \cL(X)$ containing $x$, such that the 
natural map 
$$
R_u(P) \times S \to \cL(X)
$$ 
is an open immersion. Moreover, $S$ is stable under $M := P \cap L_x$,
a Levi subgroup of both $P$ and $L_x$. It follows that $\tau^{-1}(S)$
is a locally closed $M$-stable subvariety of $X$, and the natural map
$$
R_u(P) \times \tau^{-1}(S) \to X
$$ 
is an open immersion as well. 

Together with (i), this implies that $G \cdot \tau^{-1}(S)$ is an open
subset of $X$, containing all the closed $G$-orbits. Thus,
$\tau^{-1}(S)$ meets all $G$-orbits, and it suffices to check the
desired assertion for the fibers of the restriction 
$\tau^{-1}(S) \to S$. But one shows as in the proof of Lemma
\ref{lem:reg} that $\tau^{-1}(S)$ is fixed pointwise by the derived
subgroup of $M$, and is a toric variety under a quotient of
$C(M)^0$. It follows that the irreducible components of every fiber
$\tau^{-1}(s)$ are toric varieties under a quotient of the isotropy
group $C(M)^0_s$.
\end{proof}

\begin{remarks}\label{rem:sol}
(i) The Tits map is constant if and only if $X$ is log parallelizable. 

\smallskip

\noindent
(ii) Likewise, the flag variety $\cF(X)$ is a point if and only if
every closed $G$-orbit in $X$ is an abelian variety. This is also
equivalent to the solvability of the group $G$ (or of $G_{\aff}$ by
Lemma \ref{lem:G}).

Indeed, if $G_{\aff}$ is solvable, then every closed $G_{\aff}$-orbit
in $Y$ is a point. Thus, every closed $G$-orbit in $X$ is isomorphic
to $\cA(X)$. Conversely, if $\cF(X)$ is a point, then $Y_L$ contains a
fixed point of $L$. By Proposition \ref{prop:crit}, it follows that
$L = T$ is a torus, and hence $G_{\aff} = R_u(G_{\aff}) T$ is solvable.
\end{remarks}

Next, we consider the product morphism
$\alpha \times \tau : X \to \cA(X) \times \cL(X)$
for which we obtain a generalization of the Borel--Remmert theorem: 

\begin{theorem}\label{thm:atm}
With the above notation, $\alpha \times \tau$ is surjective, and
all irreducible components of its fibers are nonsingular toric
varieties. 

Moreover, $\alpha \times \tau$ maps isomorphically every closed
$G$-orbit in $X$ to $\cA(X) \times \cF(X)$, the unique closed
$G$-orbit in $\cA(X) \times \cL(X)$.
\end{theorem}

\begin{proof}
The surjectivity of $\alpha \times \tau$ follows from Lemma
\ref{lem:at}, and the assertion on its fibers has been established in
the proof of Proposition \ref{prop:tm}.

By that proposition, any closed $G$-orbit in 
$\cA(X) \times \cL(X)$ is contained in $\cA(X) \times \cF(X)$. But the 
latter is a unique $G$-orbit, since $C^0$ acts transitively on
$\cA(X)$ and fixes pointwise $\cL(X)$.

The closed $G$-orbits in $X$ are exactly the subsets
$G\times^I Y'$, where $Y'\subseteq Y$ is a closed orbit of $I$ or,
equivalently, of $G_{\aff}$. Moreover, the restriction of $\alpha$ to
$G\times^I Y'$ is the projection to $G/I$, and $\tau$ restricts to an
isomorphism $Y' \cong \cF(X)$ by Proposition \ref{prop:tm} again. This
implies that the restriction of $\alpha \times \tau$ to 
$G \times^I Y'$ is an isomorphism to $\cA(X) \times \cF(X)$.
\end{proof}

In fact, the connected components of fibers of $\alpha \times \tau$
(regarded as reduced subschemes of $X$) are ``stable toric varieties'' 
in the sense of \cite{Al02}, as shown by the argument of Proposition
\ref{prop:tm}. Likewise, the connected components of fibers of $\tau$
are ``stable semiabelic varieties''.

The morphism $\alpha \times \tau$ may also be used to describe the 
complete log homogeneous varieties having a nonconnected boundary, 
analogously to a result of complex geometry (see \cite[I.2.6]{HO84}
and its references): 

\begin{proposition}\label{prop:ncb}
The complete log homogeneous varieties $X$ having a nonconnected
boundary are exactly the projective line bundles over the product
$A \times F$ of an abelian variety with a flag variety, obtained as
projective completions of line bundles $L \boxtimes M$, where 
$L \to A$ is a line bundle of degree $0$, and $M \to F$ is an
arbitrary line bundle. Then the Albanese map of $X$ is the projection
to $A = \cA(X)$, and the Tits map is the projection to 
$F = \cL(X) = \cF(X)$.
\end{proposition}

\begin{proof}
Let $(X,D)$ be a homogeneous pair with $X$ complete and $D$
nonconnected. Write $X = G \times^I Y$ and $D = G \times^I E$, so that
$E$ is not connected as well. Let $U \subseteq G_{\aff}$ be a maximal
unipotent subgroup. By \cite{Ho69}, the fixed point set $Y^U$ is
connected. Hence the open $G_{\aff}$-orbit $Y \setminus E$ contains
$U$-fixed points. In particular, $R_u(G_{\aff})$ fixes $Y$
pointwise; we may thus assume that $G_{\aff}$ is reductive. Let 
$J \subseteq G_{\aff}$ be the isotropy group of 
$y \in (Y \setminus E)^U$. Since $J \supseteq U$, then 
$P := N_{G_{\aff}}(J)$ is a parabolic subgroup of $G_{\aff}$ and the
quotient $P/J$ is a torus (see e.g. \cite[Sec.~7]{Ti06}). Thus, the
normalization $W$ of the closure $\overline{P \cdot y}$ is a complete
toric variety under $P/J$. Since the natural morphism
$$
\pi : G_{\aff} \times^P W \to Y
$$
is birational and $G$-equivariant, and the boundary of any complete
toric variety of dimension $\geq 2$ is connected, we must have
$\dim(W) = 1$, i.e., $P/J \cong \bG_m$ and $W \cong \bP^1$. And since
$Y$ is nonsingular and $E$ is a divisor, it follows that $\pi$ is an
isomorphism. In other words, $Y$ is the projective line bundle over 
$F := G_{\aff}/P$ associated with the principal $\bG_m$-bundle
$G_{\aff}/J \to G_{\aff}/P$ 
or with the corresponding line bundle $M \to F$. The Tits morphism is
just the structure map $Y \to F$, so that $F = \cF(X) = \cL(X)$.

As a consequence, $\alpha \times \tau$ is a projective line bundle
over $\cA(X) \times F$, associated with a line bundle of the form 
$L \boxtimes M$ for some line bundle $L \to \cA(X) =: A$. Since $L$ is
$G$-linearized, it has degree $0$.

Conversely, any such projective bundle is a complete log homogeneous
variety with boundary being the union of the zero and infinity
sections.
\end{proof}

Returning to an arbitrary complete log homogeneous variety $X$ under a
faithful action of $G$, we show that the Stein factorization of the
Tits morphism $\tau$ is intrinsic (although $\tau$ is not):

\begin{proposition}\label{prop:sf}
With the above notation, let 
$$
\CD
X @>{\sigma}>> \cS(X) @>{\varphi}>> \cL(X)
\endCD
$$ 
be the Stein factorization of the Tits morphism. Then:

\smallskip

\noindent
{\rm (i)} $\cS(X)$ is a normal projective $G$-variety, fixed pointwise
by $C^0$, and spherical under any Levi subgroup of $G_{\aff}$.

\smallskip

\noindent
{\rm (ii)} $G$ has a unique closed orbit in $\cS(X)$; it is
isomorphic to $\cF(X)$.

\smallskip

\noindent
{\rm (iii)} $\sigma$ is the morphism associated with the section ring
of the globally generated divisor $- K_X - D$, where $D$ denotes the
boundary of $X$. In particular, $\sigma$ depends only on the pair
$(X,D)$.

\smallskip

\noindent
{\rm (iv)} $\tau$ (or, equivalently, $\sigma$) is finite if and only
if $- K_X - D$ is ample. Likewise, $\alpha \times \tau$ (or 
$\alpha \times \sigma$) is finite if and only if  
$- K_X - D$ is ample relatively to $\alpha$.
\end{proposition}

\begin{proof}
(i) follows from the construction of $\sigma$.

(ii) It suffices to show that $\varphi$ restricts to an isomorphism
over $\cF(X)$. As in the proof of the above proposition, choose a
maximal unipotent subgroup $U \subseteq G_{\aff}$. Then 
$B := N_{G_{\aff}}(U)$ is a Borel subgroup of $G_{\aff}$. The fixed
point subset $\cL(X)^B$ consists of a unique point $z$, which lies in
$\cF(X)$. Moreover, $\cL(X)^U$ is connected and $B$-stable, where $B$
acts on this set through the torus $B/U$. It follows that $\cL(X)^U$
also consists of $z$. Since $\varphi$ is finite, then
$\cS(X)^U = \varphi^{-1}(z)$ (as sets). But $\cS(X)^U$ is connected as
well, so that it consists of a unique point $z'$. Since $G_{\aff,z'}$
is a subgroup of finite index of $G_{\aff,z}$ and the latter is a
parabolic subgroup of $G_{\aff}$, then
$$
\varphi^{-1}(\cF(X)) = G_{\aff} \cdot z' \to G_{\aff} \cdot z = \cF(X)
$$
is an isomorphism.

(iii) The map $op_{X,D}$ induces a surjective map
$$
\cO_X \otimes \wedge^n \fg \to 
\wedge^n \cT_X (-\log D) = \cO_X(- K_X - D),
$$
where $n = \dim(X)$. Taking global sections, we obtain a linear map
$$
\wedge^n \fg \to \Gamma(X, \cO_X(- K_X - D))
$$
with base-point-free image. The corresponding morphism
$$
X \to \bP(\wedge^n \fg^*)
$$
may be identified with $\tau$, regarded as a morphism to the Grassmann 
variety of $n$-dimensional quotients of $\fg$. This yields our
assertion. 

(iv) follows from (iii) together with Proposition \ref{prop:tm}.
\end{proof}

In view of these results, one may ask if the Stein factorization of
$\alpha \times \tau$ is $\alpha \times \sigma$, i.e., if the fibers of
the latter morphism are connected. The answer is generally negative,
as shown by the following:

\begin{example}
Let $A$ be an abelian variety and choose a point $a \in A$ of order
$2$. Consider
$$
X := A \times^F (\bP^1 \times \bP^1),
$$
where $F$ denotes the group of order $2$ acting on $A$ via translation
by $a$, and on $\bP^1 \times \bP^1$ by exchanging the two
copies. The group
$$
G := A \times \PGL_2
$$
acts on $X$ via its action on $A \times \bP^1 \times \bP^1$
by 
$$
(x,g) \cdot (y, z_1,z_2) = (x + y, g \cdot z_1, g \cdot z_2)
$$
(which commutes with $F$). Then $X$ is a projective nonsingular
variety, log homogeneous under a faithful $G$-action; the boundary is 
the prime divisor
$$
D := A \times^F \diag(\bP^1) \cong A/F \times \bP^1.
$$
Moreover, $\cA(X) = A/F$, $I = F \times \PGL_2$, $G_{\aff} = \PGL_2$
and $Y = \bP^1 \times \bP^1$. 

One checks that the Tits morphism $\tau$ may be identified with the
natural map 
$$
A \times^F (\bP^1 \times \bP^1) \to 
(\bP^1 \times \bP^1)/F \cong \bP^2.
$$
As a consequence, $\tau$ has connected fibers (isomorphic to $A$ over
$X \setminus D$, resp. to $A/F$ over $D$). But $\alpha \times \tau$ is
identified with the natural map
$$
A \times^F (\bP^1 \times \bP^1) \to A/F \times \bP^2,
$$
a double cover ramified along $D$.
\end{example}

However, $\alpha \times \sigma$ has connected fibers when $X$ is
replaced with a suitable \'etale cover, see Lemma \ref{lem:sim}
below.

\subsection{The equivariant automorphism group}
\label{subsec:aut}

We consider the group $\Aut_G(X)$, where $X \cong G \times ^I Y$
still denotes a complete log homogeneous $G$-variety. We first obtain 
a stronger version of Lemma \ref{lem:ea}:

\begin{proposition}\label{prop:ea}
$\Aut_{G_{\aff}}(Y)$ is diagonalizable and equals $\Aut_I(Y)$. 
Moreover, $\Aut_G(X)$ is commutative and preserves every
$G$-orbit in $X$. 
\end{proposition}

\begin{proof}
By Lemma \ref{lem:aut}, the subgroup 
$\Aut_{G_{\aff}}(Y) \subseteq \Aut_L(Y)$
is diagonalizable; in particular, commutative. Moreover, 
writing $I = G_{\aff} F$ as in (\ref{eqn:if}), the group $F$ acts on
$Y$ via $G_{\aff}$-equivariant automorphisms. Thus, 
$\Aut_I(Y) = \Aut_{G_{\aff}}(Y)$.

The group $\Aut_G(X)$ is generated by $C^0$ and
$\Aut_{G_{\aff}}(Y)$, and these commutative groups commute
pairwise. Thus, $\Aut_G(X)$ is commutative (alternatively, this
follows from the exact sequence (\ref{eqn:ig}) together with the
commutativity of any extension of an abelian variety by a
diagonalizable group). Moreover, $\Aut_{G_{\aff}}(Y)$ preserves all
the $G_{\aff}$-orbits in $Y$ by Lemma \ref{lem:aut}. Thus, $\Aut_G(X)$
preserves all the $G$-orbits.
\end{proof}

Next, we describe the connected equivariant automorphism group:

\begin{proposition}\label{prop:cc}
{\rm (i)} $\Aut_G(X)^0$ is a semiabelian variety, and its Albanese
variety is $\cA(X)$.

\smallskip

\noindent
{\rm (ii)} Let $X_0 \cong G/H$ be the open $G$-orbit in $X$. Then
\begin{equation}\label{eqn:iso}
\Aut_G(X)^0 \cong \Aut_G(X_0)^0 \cong (N_G(H)/H)^0 = N_G(H)^0H/H.
\end{equation}
Moreover, $N_G(H) = N_G(H^0) = N_G(\fh)$.

\smallskip

\noindent
{\rm (iii)} $\Aut_G(X)^0$ is the connected center of $\Aut^0(X,D)$,
and acts with an open orbit on any component of a general fiber of
the Tits morphism $\tau$.
\end{proposition}

\begin{proof}
(i) follows from the above proposition together with (\ref{eqn:ig}).

(ii) We claim that the (injective) restriction map 
$$
\Aut_G^0(X) \to \Aut_G^0(G/H)
$$ 
is surjective. For this, we may reduce to the case where 
$G = G_{\aff}$ and $X = Y$ by using (\ref{eqn:ig}). Then 
any element of $\Aut_G^0(G/H)$ extends to an automorphism of $X$, by
Lemma \ref{lem:aut}. This extension must be $G$-equivariant; this
proves our claim.

It remains to show that $N_G(H) = N_G(H^0)$. For this, by the
inclusions 
$$
H^0 \subseteq H \subseteq N_G(H) \subseteq N_G(H^0),
$$
it suffices to check that $N_G(H^0)/H^0$ is commutative. By Lemma
\ref{lem:G}, we may assume again that $G = G_{\aff}$; then 
$$
N_G(H^0)/H^0 = \Aut_G(G/H^0) \subseteq \Aut_L(G/H^0)
$$
and the latter group is commutative, since the $L$-variety $G/H^0$ 
(a finite covering of $G/H$) is spherical by Lemma \ref{lem:reg}.

(iii) Clearly, $\Aut_G^0(X) = \Aut_G^0(X,D)$ contains the connected
center of $\Aut^0(X,D)$. To show the equality, it suffices to check
that both groups have the same dimension. But by (ii), 
$\dim \Aut_G^0(X) = \dim N_G(\fh) - \dim H$ equals the dimension of
the general fibers of $\tau$, and the latter dimension depends only on
$(X,D)$ by Proposition \ref{prop:sf}.

Recall that the restriction $\tau\vert_{X_0}$ is the natural map 
$G/H \to G/N_G(\fh)$. Together with (\ref{eqn:iso}), this implies that
any component of a fiber of $\tau \vert_{X_0}$ is a unique
$\Aut_G(X)^0$-orbit. 
\end{proof}

We now investigate the relation between the connected equivariant
automorphism group and the morphism $\sigma: X \to \cS(X)$: 

\begin{proposition}\label{prop:pf}
Let $X$ be a complete log homogeneous variety with boundary $D$ and
put $G := \Aut^0(X,D)$. Denote by $\cE$ the subsheaf of the tangent
sheaf $\cT_{\cS(X)}$, generated by the image of $\fg$. Then $\cE$ is
locally free, i.e., the $G$-variety $\cS(X)$ is pseudo-free. Moreover,
there is an exact sequence of $G$-linearized sheaves 
\begin{equation}\label{eqn:pf}
0 \to \cO_X \otimes \fc \to \cT_X( - \log D) \to \sigma^*\cE \to 0,
\end{equation}
where $\fc$ denotes the Lie algebra of $C$.
\end{proposition}

\begin{proof}
We claim that $\fc \cap \fg_{(x)} = 0$ for all $x \in X$. Indeed, 
$\fc \cap \fg_{(x)}$ is the Lie algebra of the group $C^0_{(x)}$,
and the latter fixes pointwise $G \cdot x$ and its normal space at $x$. 
Moreover, $C^0_{(x)}$ is a torus by Lemma \ref{lem:H} and
Proposition \ref{prop:cc}(i). Thus, $C^0_{(x)}$ is trivial.

By that claim, the image of the Tits map (regarded as a morphism to
the Grassmannian $\Gr(\fg)$ of subspaces of $\fg$) consists of
subspaces $V$ such that $V \cap \fc = 0$. Let 
$\Gr^0(\fg) \subseteq \Gr(\fg)$ 
be the corresponding open subset, and $p: \Gr^0(\fg) \to \Gr(\fg/\fc)$
the projection. Then the restriction to $\Gr^0(\fg)$ of the quotient
sheaf $\cQ_{\fg}$ on $\Gr(\fg)$ fits into an exact sequence
$$
0 \to \cO_{\Gr^0(\fg)} \otimes \fc \to \cQ_{\fg} \to p^* \cQ_{\fg/\fc}
\to 0.
$$
Pulling back to $X$ yields an exact sequence 
$$ 
0 \to \cO_X \otimes \fc \to \cT_X (- \log D) \to \tau^* \cF \to 0,
$$
where $\cF$ is a $G$-linearized locally free sheaf on $\cL(X)$
(the pull-back of the quotient sheaf). In turn, this implies the exact
sequence (\ref{eqn:pf}), where $\cE = \varphi^*\cF$ for the morphism
$\varphi: \cS(X)  \to \cL(X)$. Also, $\cE$ is generated by a subspace
of global sections, the image of the composite map
$$
\fg \to \fg/\fc \to \Gamma(\Gr(\fg/\fc),\cQ_{\fg/\fc}) \to 
\Gamma(\cS(X),\cE).
$$
The restriction of (\ref{eqn:pf}) to the open orbit 
$X_0 = G \cdot x \subseteq X$
corresponds to the exact sequence of $H$-modules
$$
0 \to \fc \to \fg/\fh \to \fg/\fn_{\fg}(\fh) \to 0,
$$
where $H = G_x$. It follows that the restriction of $\cE$ to the open
orbit $G \cdot \sigma(x) \subseteq \cS(X)$ is the tangent
sheaf. Moreover, the locally free sheaf $\cE$ is a subsheaf of
$i_*(i^*\cE)$, where $i : G \cdot \sigma(x) \to \cS(X)$ denotes the
inclusion. Thus, $\cE$ may be identified with the
$\cO_{\cS(X)}$-subsheaf of derivations of the function field of
$\cS(X)$, generated by the image of $\fg$.
\end{proof}

\begin{remark}
The exact sequence (\ref{eqn:pf}) means that the relative logarithmic
tangent sheaf $\cT_{\sigma}( - \log D)$ is isomorphic to 
$\cO_X \otimes \fc$. This may be seen as a differential version of
the fact that the general fibers of $\sigma$ are $C^0$-orbits 
(Proposition \ref{prop:cc}(iii)). 

Likewise, there is an exact sequence
$$
0 \to \cO_X \otimes \fc_{\aff} \to \cT_X (- \log D) \to
(\alpha \times \sigma)^*( \cT_{\cA(X)} \boxtimes \cE) \cong
(\cO_X \otimes \fg/\fg_{\aff}) \oplus \sigma^* \cE \to 0,
$$
where $\fc_{\aff}$ denotes the Lie algebra of $C^0_{\aff}$. In other
words, the relative logarithmic tangent sheaf of 
$\alpha \times \sigma$ is isomorphic to $\cO_X \otimes \fg_{\aff}$.

\end{remark}

Finally, we consider again the morphism $\sigma : X \to \cS(X)$, under
the assumption that $\cA(X) = \cA(C^0)$ (this may be achieved by
replacing $X$ with a finite \'etale cover, see Remark \ref{rem:sim}(ii)).

\begin{lemma}\label{lem:sim}
If the isogeny $\cA(C^0) \to \cA(X)$ is an isomorphism, then 
all fibers of $\alpha \times \sigma$ are connected. Equivalently, the
Stein factorization of $\alpha \times \tau$ is $\alpha \times \sigma$.
\end{lemma}

\begin{proof}
By assumption, both isogenies $\cA(C^0) \to \cA(G)$ and 
$\cA(G) \to \cA(X)$ are isomorphisms. Thus, $C^0 \cap I$ is
connected, and $I = G_{\aff}$. So $H \subseteq G_{\aff}$ and
$C^0 \cap H = C^0_{\aff}$. Hence the restriction $\tau \vert_{X_0}$ has
fiber at the base point
$$
N_G(\fh)/H = C^0 N_{G_{\aff}}(\fh)/H \cong 
C^0 \times^{C^0_{\aff}} N_{G_{\aff}}(\fh)/H.
$$
The components of this fiber are the orbits of the group
$$
N_G(\fh)^0 H/H \cong 
C^0 \times^{C^0_{\aff}} N_{G_{\aff}}(\fh)^0 H/H.
$$
Thus, the restriction of $\alpha$ to every such component has
irreducible fibers, isomorphic to $N_{G_{\aff}}(\fh)^0 H/H$. The
assertion follows. 
\end{proof}

\subsection{Strongly log homogeneous varieties}
\label{subsec:shp}

Consider again a complete log homogeneous $G$-variety 
$X \cong G \times^I Y$. We say that $X$ is 
{\it strongly log homogeneous} if $Y$ is regular under some connected 
reductive subgroup of $G_{\aff}$.

This condition holds, for example, if $G$ is solvable (since $Y$ is
then a toric variety under a maximal torus of $G$). But it fails,
e.g., when $X = Y$ is the projective completion of the tangent bundle
$T_{\bP^2}$, and $D$ is the divisor at infinity. Then one checks
that the pair $(X,D)$ is homogeneous, and a Levi subgroup of
$\Aut^0(X,D)$ is $L := \PGL_3 \times \bG_m$, where $\PGL_3$ acts
naturally on $X$, and $\bG_m$ acts by multiplication on the fibers of
$T_{\bP^2}$. Moreover, the $L$-orbits in $X$ are the zero section
$\bP^2$, the projectivization $\bP(T_{\bP^2})$ (isomorphic to the flag
variety of $\bP^2$), and the complement $T_{\bP^2} \setminus \bP^2$;
each of them is a unique $\PGL_3$-orbit. Thus, $X$ is not regular under
$L$ or $\PGL_3$. Moreover, both closed $L$-orbits are fixed pointwise
by $\bG_m$, and they are not simultaneously regular under any proper
reductive subgroup of $\PGL_3$. So $X$ cannot be strongly log
homogeneous with boundary $D$. 

We shall obtain a criterion for strong log homogeneity which is much
simpler than our characterization of log homogeneity (Theorem
\ref{thm:str}). To motivate the following statement, observe that the
boundary of the strongly log homogeneous $G$-variety $X$ is 
$G \times^I E$, where $E$ is a subdivisor of the boundary of the
$L$-variety $Y$.

\begin{proposition}\label{prop:slh}
Let $L$ be a connected reductive group, $Y$ a complete regular
$L$-variety, $\partial Y$ its boundary, and $E \subset \partial Y$ a
subdivisor. Write
$$
\partial Y = E + D_1 + \cdots + D_r\; ,
$$
where $D_1,\ldots,D_r$ are $L$-stable prime divisors. 
Then the pair $(Y,E)$ is (strongly) log homogeneous if and only if 
$D_1, \ldots, D_r$ are all generated by their global sections.  

Under one of these assumptions, we have an exact sequence 
\begin{equation}\label{eqn:httn}
0 \to \Gamma(Y, \cT_Y( - \log \partial Y)) \to 
\Gamma(Y, \cT_Y( - \log E)) \to 
\bigoplus_{i=1}^r \Gamma(D_i, \cN_{D_i/Y}) \to 0,
\end{equation}
where $\cN_{D_i/Y}$ denotes the normal sheaf. Moreover,
$H^j(Y, \cT_Y( - \log E)) = 0$ for all $j \geq 1$.
\end{proposition}

\begin{proof}
We adapt arguments from \cite[Sec.~4.1]{BB96}; we provide details
for  completeness. Consider the map
$$
p: \cT_Y( - \log E) \to \bigoplus_{i=1}^r \cN_{D_i/Y}
$$
obtained by taking the direct sum of the natural maps 
$\cT_Y \to \cN_{D_i/Y}$ and then restricting to $\cT_Y (- \log E)$.
Clearly, the kernel of $p$ is $\cT_Y( - \log \partial Y)$. 
Moreover, by using a local system of parameters, one checks that $p$
is surjective. This yields a short exact sequence
\begin{equation}\label{eqn:ttn}
0 \to \cT_Y( - \log \partial Y) \to \cT_Y( - \log E) \to 
\bigoplus_{i=1}^r \cN_{D_i/Y} \to 0.
\end{equation}
Taking the associated long exact sequence of cohomology groups and
using the vanishing of $H^j(Y, \cT_Y( - \log \partial Y))$ for 
$j \ge 1$ (a consequence of \cite[Thm.~4.1]{Kn94b}), we obtain the
short exact sequence (\ref{eqn:httn}) and isomorphisms
\begin{equation}\label{eqn:hh}
H^j(Y, \cT_Y( - \log E)) \cong \bigoplus_{i=1}^r 
H^j(D_i, \cN_{D_i/Y})
\end{equation}
for all $j \geq 1$. 

From the exact sequences (\ref{eqn:httn}), (\ref{eqn:ttn}) and the
global generation of $\cT_Y (- \log \partial Y)$ (Corollary
\ref{cor:rlh}), it follows that $\cT_Y(- \log E)$ is globally
generated if and only if so is each $\cN_{D_i/Y}$. Moreover, from the
standard exact sequence 
$$
0 \to \cO_Y \to \cO_Y (D_i) \to \cN_{D_i/Y} \to 0
$$
and the vanishing of $H^j(Y, \cO_Y)$ for $j \geq 1$ (which follows
e.g. from \cite[Cor.~31.1]{Ti06}), we obtain an exact sequence 
$$
0 \to k \to \Gamma(Y, \cO_Y (D_i)) \to \Gamma(D_i, \cN_{D_i/Y}) \to 0
$$
and isomorphisms
\begin{equation}\label{eqn:cong}
H^j(Y, \cO_Y (D_i)) \cong H^j(D_i, \cN_{D_i/Y})
\end{equation}
for all $j \geq 1$. Hence the global generation of
$\cN_{D_i/Y}$ is equivalent to that of $\cO_Y(D_i)$. The latter
implies the vanishing of $H^j(Y, \cO_Y(D_i))$ for all $j \ge 1$, since
any globally generated invertible sheaf on a complete spherical
variety has vanishing higher cohomology (see \cite[Cor.~31.1]{Ti06}
again). By (\ref{eqn:hh}) and (\ref{eqn:cong}), this implies in turn
the vanishing of the $H^j(Y, \cT_Y(- \log E))$. Also, we have shown that
the global generation of $\cT_Y(- \log E)$ is equivalent to that of
$D_1, \ldots, D_r$.
\end{proof}

This criterion may be formulated in combinatorial terms by using
the characterization of globally generated divisors on spherical
varieties, see e.g. \cite[Sec.~17]{Ti06}. For example, if $L$ is a
torus, i.e., $Y$ is a toric variety, let $\Delta$ be the corresponding 
fan, and $\rho_1, \ldots, \rho_r \in \Delta$ the rays associated with
$D_1,\ldots,D_r$. Then the global generation of $D_i$ is equivalent to
the convexity of the union of all the cones of $\Delta$ that do not
contain $\rho_i$ (in other words, the complement of the star of
$\rho_i$).

Proposition \ref{prop:slh} is the main ingredient in the classification
of complete log homogeneous surfaces (see Remark \ref{rem:lh}(iii) for
the much easier classification of log homogeneous curves):

\begin{proposition}\label{prop:cs}
Up to isomorphism, the pairs $(X,D)$, where $X$ is a complete log
homogeneous surface with boundary $D$, are those in the following
list:

\smallskip

\noindent
{\rm (parallelizable)} $(A,\emptyset)$, where $A$ is an abelian
surface.

\smallskip

\noindent
{\rm (homogeneous)} $(E \times \bP^1,\emptyset)$, where $E$ is an
elliptic curve, $(\bP^1 \times \bP^1,\emptyset)$, and
$(\bP^2,\emptyset)$. 

\smallskip

\noindent
{\rm (log parallelizable)} $X$ is the projective completion of a line
bundle of degree $0$ on an elliptic curve, and $D$ is the union of the
zero and infinity sections.

\smallskip

\noindent
{\rm (log homogeneous)}

{\rm (a)} $(E \times \bP^1, E \times \{\infty\})$, where $E$ is an
elliptic curve.

{\rm (b)} $(\bP^1 \times \bP^1, \diag(\bP^1))$.

{\rm (c)} $(\bP^2,C)$, where $C$ is a conic.

{\rm (d)} $(\bP^2,D)$, where $D$ is a subdivisor of the boundary of the
toric variety $\bP^2$, i.e., a union of coordinate lines.

{\rm (e)} $(\bP^1 \times \bP^1, D)$, where $D$ is a subdivisor of the
boundary of the toric variety $\bP^1 \times \bP^1$. 

{\rm (f)} $(\bF_n,D)$, where $\bF_n$ is the rational ruled surface of
index $n \geq 1$ and $D$ is a subdivisor of the boundary of the toric 
variety $\bF_n$ that contains the unique curve $C_{-n}$ of
self-intersection $-n$.

{\rm (g)} All pairs obtained from those in {\rm (d)}, {\rm (e)}, 
{\rm (f)} by successively blowing up intersection points of boundary
curves.
\end{proposition}

We refer to \cite[Sec.~1.7]{Od88} for a description of $\bF_n$
regarded as a toric surface. We outline the proof of Proposition
\ref{prop:cs}, leaving the easy but rather long details to the
reader. The parallelizable resp.~homogeneous cases are
straightforward, and the log parallelizable case follows either from
Theorem \ref{thm:lpt} or from Proposition \ref{prop:ncb}. So we may
assume that $D$ is nonempty and connected. If $D$ is irreducible, then
by looking at the minimal model of $X$, one obtains Cases (a), (b),
(c), (d) and (e) (where $D$ is a line), as well as (f) 
(where $D = C_{-n}$). Otherwise, $X$ contains a fixed point of $G$,
namely, the intersection of two boundary curves. By Theorem
\ref{thm:atm} and Remark \ref{rem:sol}(ii), it follows that $G$ is
affine and solvable. In particular, $X$ is a toric surface and $D$ is
a subdivisor of the toric boundary. Then the statement follows from
Proposition \ref{prop:slh} by elementary arguments of two-dimensional
toric geometry.

\begin{remarks}\label{rem:fin}
(i) In Cases (a)-(e), the only possibility for the acting group $G$ is
the full automorphism group $\Aut^0(X,D)$. But this does not
extend to Case (f): consider, for example, the pair 
\begin{equation}\label{eqn:bad}
(X,D) = (\bF_n, C_{-n} + F_0 + F_{\infty}),
\end{equation}
where $F_0$ and $F_{\infty}$ are the fibers of the ruling 
$$
\bF_n = \bP( \cO_{\bP^1} \oplus \cO_{\bP^1}(n)) \to \bP^1
$$ 
at the corresponding points. Then $\Aut^0(X,D)$ is the semidirect
product of the group $U := H^0(\cO_{\bP^1}(n))$ (acting by
translations) with a torus $T$ of dimension $2$ (acting on $U$ via the
image of a maximal torus of $\GL_2$). Moreover, we may take for $G$
the semidirect product of $V$ with $T$, where $V \subseteq U$ is any
$T$-stable subspace without base points; there are $2^{n-1}$ such
subspaces.

\smallskip

\noindent
(ii) For $(X,D)$ as in (\ref{eqn:bad}) and $G = \Aut^0(X,D)$, the
image of $\sigma$ is the normal surface obtained from $\bF_n$ by
contracting $C_{-n}$. Indeed, the group $\Aut_G(X)^0$ is trivial, so
that $\sigma$ is birational. Moreover, $\sigma$ maps the two $G$-fixed
points $F_0 \cap C_{-n}$ and $F_{\infty}\cap C_{-n}$ to the same point. 

In particular, the variety $\cS(X)$ is singular if $n \geq 2$. One may
also check that $\varphi : \cS(X) \to \cL(X)$ is bijective for all
$n$, so that $\cL(X)$ is also singular if $n \geq 2$.
\end{remarks}

\end{document}